\documentclass{amsart}

\theoremstyle{definition}

\def\half{{\textstyle{\frac12}}}
\def\third{{\textstyle{\frac13}}}
\def\fourth{{\textstyle{\frac14}}}

\def\u{\mathsf u} 
\def\v{\mathsf v} 
\def\w{\mathsf w} 
\def\y{\mathsf y} 

\def\^{$\hat{\ }$} 

\def\={\!\!\!\!=\!\!\!\!}

\usepackage{amssymb} 

\title[Finsler metrics of constant positive curvature]
      {FINSLER METRICS OF CONSTANT POSITIVE \\ 
       CURVATURE ON THE LIE GROUP $S^3$}

\author[Bao and Shen]{DAVID BAO and Z. SHEN}


\address{Department of Mathematics, University of Houston \\
         Houston, TX 77204-3476, USA}
\address{Department of Mathematical Sciences, IUPUI \\
         Indianapolis, IN 46202-3216, USA}

\email{bao@math.uh.edu}
\email{zshen@math.iupui.edu}

\begin{document}

\begin{abstract}
Guided by the Hopf fibration, we single out a family 
(indexed by a positive constant $K$) of right invariant 
Riemannian metrics on the Lie group $S^3$. Using the Yasuda--Shimada theorem 
as an inspiration, we determine for each $K > 1$ a privileged right invariant 
Killing field of constant length. Each such Riemannian metric pairs with the 
corresponding Killing field to produce a $y$-global and {\it explicit} Randers 
metric on $S^3$. Using the machinery of spray curvature and Berwald's 
formula for it, we prove directly that the said Randers metric has constant 
positive flag curvature $K$, as predicted by the Yasuda--Shimada theorem. 
We also explain why this family of Finslerian space forms is {\it not} 
projectively flat.  
\end{abstract}

\maketitle

\section{Introduction}

A Finsler metric $F$ is a family of `norms' on a manifold $M$, one on each 
tangent space $T_xM$. These `norms' are typically only positively homogeneous 
of degree one, whereas the norms used in functional analysis are absolutely 
homogeneous. There are also the usual smoothness and strong convexity 
assumptions (see for instance \cite{BCS}) on the slit tangent bundle 
$TM \smallsetminus 0$. In a large number of examples, especially the ones of 
physical origin, these technical requirements are only satisfied on open 
cones in $TM \smallsetminus 0$. If the Finsler metric is smooth and strongly 
convex on the {\it entire} slit tangent bundle $TM \smallsetminus 0$, it is 
said to be $y$-global. Examples that are both $y$-global and geometrically 
significant are highly sought after. 

In Riemannian geometry, one has the concept of sectional curvature. Its 
analogue in Finsler geometry is called the flag curvature. Flag curvatures 
are more easily accessed in some settings through spray curvatures. 
These objects have now been given detailed treatments in textbooks and 
monographs such as \cite{R}, \cite{BCS}, \cite{AIM}, and \cite{S1}. Finsler 
spaces of constant flag curvature are, just like their Riemannian 
counterparts, known as space forms. However, unlike Riemannian geometry, 
Finslerian space forms do not just arise from three standard models. 
This phenomenon is related to the fact that on $\mathbb R^n$, there are many 
isometry classes of norms, but only one isometry class of inner products. 
 
Akbar-Zadeh \cite{AZ} showed that if the metrics in question are 
geodesically complete and the growth of their Cartan tensors 
are suitably constrained, then flat and negatively curved space forms 
are fairly well understood. See \cite{BCS} for a leisurely exposition. 
Relax any of those two hypotheses and one encounters intriguing spaces, 
like for example the Finslerian Poincar\'e disc discussed in \cite{BCS}. 

This paper deals with $y$-global Finsler spaces of constant positive flag 
curvature, and puts us in an even more esoteric landscape. 
In this regard, we have several pioneering works of Bryant's. Two    
dimensional examples are treated in \cite{Br1}, \cite{Br2}, \cite{Br3};  
higher dimensional ones are discussed in his lectures. As far as we know, all 
of Bryant's examples are projectively flat (meaning that their geodesics are 
straight lines in certain coordinate systems), and none of them is of the 
Randers vintage. 

Randers spaces are Finsler spaces constructed from just two pieces of 
familiar data: a Riemannian metric and a differential $1$-form, both 
globally defined on an underlying smooth manifold. As such, they are 
possibly the best stepping stones from the Riemannian realm to the Finslerian 
territory. 
\begin{quote} 
The goal of this paper is to produce, for each constant 
$K > 1$, an explicit example of a compact boundaryless (non-Riemannian) 
Randers space that has constant positive flag curvature $K$, 
and which is not projectively flat. 
\end{quote} 

For that purpose, we turn to the Yasuda--Shimada theorem \cite{YS}.  
This result was published in the late 70s, and classifies Randers metrics 
of constant positive and constant negative flag curvature. Incidentally, 
flat Randers metrics are necessarily locally Minkowskian. For the positively 
curved case, the Yasuda--Shimada theorem gives four mathematical criteria on 
the Riemannian metric and the $1$-form. These criteria are made precise in 
\S 2. 
\begin{quote} 
For each $K > 1$, we select a manifold $M$ judiciously, then solve those four 
criteria for a Riemannian metric $\tilde a$ and a nonzero 
$1$-form $\tilde b$, both living on $M$. This is carried out in \S 3--\S 5. 
By doing so, we have shown that the Yasuda--Shimada criteria are non-vacuous. 
The simplest choice of $M$ turns out to be the Lie group $S^3$. Similar 
promise \cite{Ro} holds for $S^{2n+1}$, mainly because every odd-dimensional 
sphere admits the Hopf fibration. 
\end{quote} 

After the construction, we could in principle appeal to the Yasuda--Shimada 
theorem to conclude that our Randers spaces have constant positive flag 
curvature. However, the two published proofs (\cite{YS} and \cite{M}) of the 
theorem in question both contain arguments that we find difficult to follow. 
In order to render our paper logically self-contained, we prove directly that 
each said example indeed has constant positive flag curvature $K$. The method 
we use is quite different from the ones in \cite{YS} and \cite{M}. It is 
geometrical and brings out the spectacular power behind a formula of 
Berwald's. It also reaffirms our belief that a happy synergy awaits Finsler 
geometry and modern computing. 

To that end, we first review the concept of spray curvatures, and its 
close relationship with the flag curvature, in \S 6. We then give Berwald's 
formula for the spray curvature. This formula has proved to be 
exceptionally useful in computational Finsler geometry, whether one carries 
out the computations on machine or by hand. 

In the Appendix, we write a Maple program based on Berwald's formula and use 
it to generate numerical data, supporting the contention that our Randers 
metrics have constant positive flag curvature $K$. This type of experimental 
evidence provides the faith which sustains our direct {\it proof} (by hand) 
in \S 7. The hand computations in this direct proof are facilitated by a 
refined understanding \cite{S1} of the spray curvature, when the latter is 
adapted to the Randers setting. 

We believe that this refined understanding also holds the key to a distinctly 
new proof of the Yasuda--Shimada theorem. This theorem has been a reliable 
source of beautiful examples. For instance, the Finslerian Poincar\'e disc we 
cited above has a Yasuda--Shimada pedigree as well. Given that, it is 
desirable to have a more accessible and {\it geometric} proof of the theorem. 
The outcome of such an endeavor will be reported elsewhere.  

Finally, in \S 8, we give a rather preliminary discussion on the geodesics 
of our Randers metrics. Specifically, we invoke standard results about 
the projective Weyl and Douglas tensors to deduce that there is no coordinate 
system in which the geodesics appear as straight lines.  
Thus, as promised, our Randers spaces are not projectively flat.

\bigskip

\section{Randers spaces and the Yasuda--Shimada theorem} 

Randers metrics were introduced by Randers in 1941 \cite{Ra} in the context 
of general relativity. They play a prominent role in Ingarden's study of 
electron optics (see his treatment of the subject in \cite{AIM}). 
Mathematically, in spite of the wide range of non-Riemannian phenomena they 
are capable of producing, Randers spaces are Finsler spaces built from data 
that are quite familiar to all differential geometers:  
\begin{itemize} 
\item a Riemannian metric 
      $\tilde a := \tilde a_{ij} \, dx^i \otimes dx^j$ 
      on a smooth $n$-dimensional manifold $M$, and 
\smallskip 
\item a differential $1$-form $\tilde b := \tilde b_i \, dx^i$ on $M$. 
      This $\tilde b$ is sometimes called a drift $1$-form. 
\end{itemize}  
Together these objects define a Finsler metric $F$ in a simple way: 
\begin{equation*}
 F(x, y) := \alpha (x, y) + \beta (x, y) ,
\end{equation*}  
where 
\begin{equation*}
 \begin{aligned}  
  \alpha (x, y) &:= \sqrt{ \, \tilde a_{ij \, (x)} \, y^i y^j \, } \\ 
  \beta (x, y)  &:= \tilde b_{i \, (x)} \, y^i . 
 \end{aligned}   
\end{equation*}
Here, $x$ stands for points on the manifold $M$, and $y \in T_xM$ denote 
tangent vectors based at $x$. Those tangent space coordinates $y^i$ typically 
come from the expansion $y = y^i \frac{\partial}{\partial x^i}$ in terms 
of a local coordinate basis. They can also arise, as $\y^p$, from the 
expansion $y = \y^p e_p$ in terms of an $\tilde a$-orthonormal frame 
field on $M$. Both scenarios are exemplified in \S 5. 

Due to the presence of $\beta$, the Randers metric $F := \alpha + \beta$ is 
generally only positively homogeneous of degree one in $y$: 
$F(x, cy) = c \, F(x, y)$ for all positive $c$. A Randers 
metric cannot be absolutely homogeneous [$F(x, cy) = |c| F(x, y)$] unless 
$\tilde b = 0$, in which case $F$ is Riemannian. 

The fundamental tensor is formally analogous to the metric tensor in 
Riemannian geometry. It is defined as 
\begin{equation*} 
 g_{ij} := \half \left( F^2 \right)_{y^i y^j} \, ,
\end{equation*} 
where we have used $y^i$, $y^j$ as subscripts to signify 
partial differentiation. 
Almost by inspection, we have 
\begin{equation*}
 \begin{aligned}   
  \tilde \ell_i 
  &:=  
  \alpha_{y^i} = \frac{\tilde a_{ij} \, y^j}{\alpha} , \\   
  \ell_i 
  &:=  
  F_{y^i} = \tilde \ell_i + \tilde b_i . 
 \end{aligned} 
\end{equation*}  
The fundamental tensor can then be expressed as 
\begin{equation*} 
 g_{ij} 
 =  
 \frac{F}{\alpha} 
 \left( \tilde a_{ij} - \tilde \ell_i \, \tilde \ell_j \right) 
 +  
 \ell_i \, \ell_j . 
\end{equation*}   
 
Since $\beta (x, y)$ is linear in $y$, it cannot possibly have a fixed sign. 
The size of $\tilde b$ therefore needs to be controlled if $F$ is to be 
positive on $TM \smallsetminus 0$. We also want the fundamental tensor to be 
positive definite. It turns out that both these properties hold if and only if 
\begin{equation*} 
 \| \, \tilde b \, \| := \sqrt{ \, \tilde b_i \, \tilde b^i \, } < 1 , 
\end{equation*}  
where 
\begin{equation*} 
\tilde b^i := \tilde a^{ij} \, \tilde b_j .
\end{equation*}  
So, the drift $1$-form $\tilde b$ of Randers spaces must be required to have 
Riemannian norm strictly smaller than $1$ everywhere. See \cite{BCS} or 
\cite{AIM}.    

The notion of flag curvature makes sense for all Finsler spaces. It is 
constructed from the $hh$-curvature tensor of one's favorite Finsler 
connection. Happily, the resulting flag curvature is independent of which 
standard connection one is using, be it Berwald's, Cartan's, Chern's, 
or Hashiguchi's, just to name a few. Flag and spray curvatures will be 
reviewed in \S 6 rather than here, in order to effect a more streamlined 
exposition. 

The Yasuda--Shimada theorem \cite{YS} says that for a Randers metric to 
have constant {\it positive} flag curvature $K$, the following four criteria 
are both necessary and sufficient: 

\begin{itemize} 
\item The $1$-form $\tilde b$ is a Killing field of the 
      Riemannian metric $\tilde a$. 
      \begin{equation*}
       \tilde b_{i|j} + \tilde b_{j|i} = 0 . 
      \end{equation*}
\item The Riemannian norm of this Killing field must be constant, 
      besides being globally less than $1$.  
      \begin{equation*} 
       1 > \| \, \tilde b \, \| ^2 
       :=  
       \tilde b_i \, \tilde b^i \ \ \text{is constant}. 
      \end{equation*}
\item Second order covariant derivatives of $\tilde b$ are to have the 
      specific form 
      \begin{equation*} 
       \tilde b_{i|j|k} 
       =  
       K \left( \tilde a_{ik} \, \tilde b_j 
              - \tilde a_{jk} \, \tilde b_i \right) . 
      \end{equation*}
\item The Riemann curvature tensor $\tilde R_{hijk}$ of the metric $\tilde a$ 
      must be given by the special formula  
      $$\begin{aligned}  
        \ \ \ \ 
        &- K \left( 1 - \| \, \tilde b \, \|^2 \right) \, 
             \tilde a_{hj} \, \tilde a_{ik} 
         - K \left( \tilde a_{hj} \, \tilde b_i \, \tilde b_k 
                    +  
                    \tilde a_{ik} \, \tilde b_h \, \tilde b_j \right) 
         + \tilde b_{h|j} \, \tilde b_{i|k}                          \\ 
        &+ K \left( 1 - \| \, \tilde b \, \|^2 \right) \, 
             \tilde a_{hk} \, \tilde a_{ij} 
         + K \left( \tilde a_{hk} \, \tilde b_i \, \tilde b_j 
                    +  
                    \tilde a_{ij} \, \tilde b_h \, \tilde b_k \right) 
         - \tilde b_{h|k} \, \tilde b_{i|j}                          \\ 
        &+ 2 \, \tilde b_{h|i} \, \tilde b_{j|k} . 
       \end{aligned}$$ 
\end{itemize} 
Here, the vertical slash $({\cdots})_{i|j}$ denotes covariant differentiation 
on $M$, taken with respect to the Levi-Civita (Christoffel) connection of the 
Riemannian metric $\tilde a$. 

In order to obtain a non-Riemannian Randers metric from these four criteria, 
the drift $1$-form needs to be nowhere zero because it has constant length. 
Limiting our search to compact oriented manifolds $M$ without boundary, 
we deduce from the Poincar\'e--Hopf index theorem that the Euler 
characteristic $\chi (M)$ must vanish.   

In two dimensions, the only candidate for $M$ is therefore the torus. On  
this manifold, the most easily visualized Riemannian metric $\tilde a$ is the 
one that gives the torus of revolution in $\mathbb R^3$. In that case, a 
straightforward calculation shows that every constant length Killing field 
$\tilde b$ is identically zero. So there is no non-Riemannian Randers metric 
with constant positive flag curvature on the torus of revolution.   

It is not obvious whether the above conclusion again results if we use other 
Riemannian metrics $\tilde a$ on the torus. If we are in the boundaryless 
category, there are telltale indications 
(\cite{S1}, \cite{S2}) that complete, positively homogeneous [as opposed to 
the more restrictive absolute homogeneity $F(x, cy) = |c| F(x, y)$]  
Finsler metrics of constant positive flag curvature can only be supported on 
manifolds homeomorphic to spheres.

\bigskip

\section{The Lie group $S^3$ and its Hopf fibration} 

Our goal is to find a non-Riemannian Randers metric by `solving' the four 
criteria stated in the Yasuda--Shimada theorem. This means we must produce a 
manifold $M$ and, living globally on it, a special Riemannian metric 
$\tilde a$ and a nonzero differential $1$-form $\tilde b$ that is 
appropriately coupled to $\tilde a$. As explained near the end of \S 2, the 
base manifold $M$ must have zero Euler characteristic. Working in odd 
dimensions and staying in the boundaryless category automatically satisfies 
this topological constraint. 

The simplest compact boundaryless oriented $3$-manifold is the Lie group 
$S^3$, which is also a circle bundle over $S^2$. We briefly discuss these 
properties for the sole purpose of setting some notation. 
\begin{itemize} 
\item As a Riemannian manifold, $S^3$ is the standard unit sphere in Euclidean 
      $\mathbb R^4$, whose points $x$ have Cartesian coordinates 
      $(x^o, x^1, x^2, x^3)$. In other words, $S^3$ can be characterized by 
      \begin{equation*} 
       |x| := \sqrt{ \, (x^o)^2 + (x^1)^2 + (x^2)^2 + (x^3)^2 \, } = 1 . 
      \end{equation*}  
\item As a Lie group, $S^3$ consists of the unit quaternions 
      \begin{equation*} 
       x = x^o \textbf{1} + x^1 \textbf{i} + x^2 \textbf{j} + x^3 \textbf{k}  
      \end{equation*} 
      in the space $\mathbb H$ of quaternions. Its group structure is 
      isomorphic to $SU(2)$. Multiplication between 
      quaternions is non-commutative and is governed by the formal rules 
      \[
      \begin{array}{lcllcllcl}   
         \textbf{i} \textbf{i} & \= & -\textbf{1} , \ \  & 
         \textbf{j} \textbf{j} & \= & -\textbf{1} , \ \  &  
         \textbf{k} \textbf{k} & \= & -\textbf{1} ,             \\  
         \textbf{i} \textbf{j} & \= & +\textbf{k} , \ \  &  
         \textbf{j} \textbf{k} & \= & +\textbf{i} , \ \  &  
         \textbf{k} \textbf{i} & \= & +\textbf{j} ,             \\  
         \textbf{j} \textbf{i} & \= & -\textbf{i} \textbf{j}, \ \ & 
         \textbf{k} \textbf{j} & \= & -\textbf{j} \textbf{k}, \ \ &
         \textbf{i} \textbf{k} & \= & -\textbf{k} \textbf{i},      
      \end{array}
      \] 
      where the object \textbf{1} acts just like the real number $1$. In view 
      of the remarkable formula $|x \breve x| = |x| \, |\breve x|$, 
      the unit quaternions are closed under the defined multiplication. Also, 
      $x \bar x = |x|^2 = \bar x x$, with  
      \begin{equation*} 
       \bar x := x^o \textbf{1} - x^1 \textbf{i} 
                                - x^2 \textbf{j} 
                                - x^3 \textbf{k} . 
      \end{equation*} 
      So the group inverse of any $x \in S^3$ is its conjugate: 
      \begin{equation*} 
       x^{-1} = \bar x . 
      \end{equation*}
\end{itemize}

As a circle bundle, $S^3$ has the celebrated Hopf fibration. For that 
purpose, it is best described as the set of 
\begin{equation*} 
 (z,w) := ( \, x^o + i x^1 \, , \, x^2 + i x^3 \, ) 
\end{equation*} 
in $\mathbb C \times \mathbb C$ satisfying $|z|^2 + |w|^2 = 1$. 
Incidentally, with this notation, the aforementioned isomorphism between 
$S^3$ and $SU(2)$ reads
$$
(z, w) 
\leftrightarrow 
\begin{pmatrix}
      z  &      w   \\ 
- \bar w & \bar z  
\end{pmatrix} \ . 
$$
The unit complex numbers $e^{i \theta} \in S^1$ act on $S^3$ by  
$$(z,w) 
  \mapsto 
  \left( \, e^{i \theta} z \, , \, e^{i \theta} w \, \right) . $$ 
The quotient manifold is the same as $\mathbb C^2$ mod first the length and 
then the argument of an arbitrary complex number. So it is $\mathbb CP^1$ 
which, by a suitable stereographic projection (see \cite{F}), is precisely 
the complex manifold $S^2$. Thus $S^3$ is a circle bundle over $S^2$. This 
bundle is non-trivial because $S^3$ is simply connected, whereas 
$S^2 \times S^1$ is not. 
  
Using a curve of angles $\theta(t)$ with $\theta(0) = 0$, one calculates the 
infinitesimal generator 
$$\frac{d}{dt}_{| t=0} 
  \left( \, e^{i \theta(t)} z \, , \, e^{i \theta(t)} w \, \right)$$ 
of the said $S^1$ action. The answer is $\theta^\prime(0)$ times the 
following distinguished vector field 
\begin{equation*} 
 E_1 := - x^1 \partial_o + x^o \partial_1 
        - x^3 \partial_2 + x^2 \partial_3 ,  
\end{equation*} 
where $\partial_\mu$ abbreviates $\partial_{x^\mu}$. This $E_1$ is 
therefore tangent to the $S^1$ fibres. It is complemented by  
\begin{equation*} 
 \begin{aligned} 
  E_2 &:= - x^2 \partial_o + x^3 \partial_1 
          + x^o \partial_2 - x^1 \partial_3     \\  
  E_3 &:= - x^3 \partial_o - x^2 \partial_1 
          + x^1 \partial_2 + x^o \partial_3   
 \end{aligned} 
\end{equation*} 
to give a global orthonormal frame field on $S^3$ with 
$$[E_1, E_2] = -2 E_3 , \ \ 
  [E_2, E_3] = -2 E_1 , \ \ 
  [E_3, E_1] = -2 E_2 . $$
At any point $x \in S^3$, the values of $E_1$, $E_2$, $E_3$ are 
respectively equal to the quaternion products $\textbf{i} x$, $\textbf{j} x$, 
and $\textbf{k} x$. This fact can be used to check that our frame field is 
right invariant. 
 
The natural dual of $\{ E_1, E_2, E_3 \}$ is the following 
globally defined right invariant orthonormal coframe on $S^3$: 
\begin{equation*} 
 \begin{aligned} 
  \Theta^1 &:= - x^1 dx^o + x^o dx^1 - x^3 dx^2 + x^2 dx^3 ,   \\ 
  \Theta^2 &:= - x^2 dx^o + x^3 dx^1 + x^o dx^2 - x^1 dx^3 ,   \\
  \Theta^3 &:= - x^3 dx^o - x^2 dx^1 + x^1 dx^2 + x^o dx^3 .   
 \end{aligned} 
\end{equation*} 
The standard metric on $S^3$ is thus  
$\Theta^1 \otimes \Theta^1 + 
 \Theta^2 \otimes \Theta^2 + 
 \Theta^3 \otimes \Theta^3$. 
Also, one has 
$$d \Theta^1 = 2 \, \Theta^2 \wedge \Theta^3 , \ \ \ 
  d \Theta^2 = 2 \, \Theta^3 \wedge \Theta^1 , \ \ \ 
  d \Theta^3 = 2 \, \Theta^1 \wedge \Theta^2 . $$ 

Motivated by the treatment of the Hopf fibration in \cite{GLP}, we consider 
the following Riemannian metric on $S^3$: 
\begin{equation*} 
 \tilde a := \epsilon^2 \, \Theta^1 \otimes \Theta^1 + 
                           \Theta^2 \otimes \Theta^2 + 
                           \Theta^3 \otimes \Theta^3 , 
\end{equation*} 
where $\epsilon$ is a positive constant (typically different from $1$). 
This $\tilde a$ modifies the standard metric on $S^3$ by introducing a 
dilation along the $S^1$ fibres. The Lie derivative of $\tilde a$ can be 
calculated by using the Cartan formula 
$\mathcal L_X \Theta = i_X d \Theta + d i_X \Theta$ on  
$\Theta^1$, $\Theta^2$, and $\Theta^3$. We find that every constant multiple 
of $E_1$ is a Killing vector field of $\tilde a$, whereas $E_2$, $E_3$ are 
not Killing fields unless $\epsilon$ happens to be $1$. 

The Riemannian metric $\tilde a$ admits its own orthonormal frame field 
$$e_1 := \frac{1}{\epsilon} \, E_1 , \ \ 
  e_2 :=                       E_2 , \ \ 
  e_3 :=                       E_3 . $$
The natural dual consists of 
$$\omega^1 := \epsilon \, \Theta^1 , \ \ 
  \omega^2 :=             \Theta^2 , \ \ 
  \omega^3 :=             \Theta^3 , $$
with 
$$d \omega^1 = 2 \epsilon         \, \omega^2 \wedge \omega^3 , \ \ \ 
  d \omega^2 = \frac{2}{\epsilon} \, \omega^3 \wedge \omega^1 , \ \ \ 
  d \omega^3 = \frac{2}{\epsilon} \, \omega^1 \wedge \omega^2 . $$ 

Relative to this $\tilde a$-orthonormal frame field, the 
Levi-Civita (Christoffel) connection of $\tilde a$ consists of a 
skew-symmetric $3 \times 3$ matrix of $1$-forms $\omega_q^{\ p}$. 
These are obtained by solving the structural equations 
$d \omega^p = \omega^q \wedge \omega_q^{\ p}$  
and $\omega_{pq} = - \omega_{qp}$. Here, $\omega_{qp}$ means 
$\omega_q^{\ s} \, \delta_{sp}$, which is 
numerically the same as $\omega_q^{\ p}$. We find that 
$$
\begin{pmatrix}
\omega_1^{\ 1} & \omega_1^{\ 2} & \omega_1^{\ 3} \\ 
\omega_2^{\ 1} & \omega_2^{\ 2} & \omega_2^{\ 3} \\ 
\omega_3^{\ 1} & \omega_3^{\ 2} & \omega_3^{\ 3} 
\end{pmatrix} 
= 
\begin{pmatrix} 
0  & - \epsilon \omega^3 & \epsilon \omega^2 \\ 
\epsilon \omega^3 & 0 & \left(\epsilon-\frac{2}{\epsilon}\right) \omega^1 \\ 
- \epsilon \omega^2 & \left(\frac{2}{\epsilon}-\epsilon\right) \omega^1 & 0 
\end{pmatrix} \ . 
$$

The curvature $2$-forms and the Riemann curvature tensor are related by 
$$d \omega_q^{\ p} - \omega_q^{\ s} \wedge \omega_s^{\ p} 
  \ = \ 
  \half \, \tilde R_{q \ rs}^{\ p} \ \omega^r \wedge \omega^s , $$ 
where our convention for indices on $\tilde R_{q \ rs}^{\ p}$ follows that in 
\cite{BCS}. Since we are in an orthonormal frame, $\tilde R_{q \ rs}^{\ p}$ is 
numerically equal to $\tilde R_{qprs}$. Computations give 
\[
\begin{array}{lcllcllcl}   
 \tilde R_{1212} & \= &     - \epsilon^2 , \ \  & 
 \tilde R_{1213} & \= &                0 , \ \  &  
 \tilde R_{1223} & \= &                0 ,             \\  
 \tilde R_{1312} & \= &                0 , \ \  & 
 \tilde R_{1313} & \= &     - \epsilon^2 , \ \  &  
 \tilde R_{1323} & \= &                0 ,             \\  
 \tilde R_{2312} & \= &                0 , \ \  & 
 \tilde R_{2313} & \= &                0 , \ \  &  
 \tilde R_{2323} & \= & 3 \epsilon^2 - 4 .              
\end{array}
\] 
All other components of the Riemann curvature of $\tilde a$ are obtained from 
these by standard properties. Namely, $\tilde R_{qpsr} = - \tilde R_{qprs}$, 
$\tilde R_{pqrs} = - \tilde R_{qprs}$, and the block symmetry 
$\tilde R_{rsqp} = \tilde R_{qprs}$.  

\bigskip

\section{A privileged Killing field on $S^3$} 

In \S 3, we discussed the Hopf fibration of $S^3$ and introduced the 
vector field $E_1$ that is tangent to the $S^1$ fibres. The natural 
dual of this $E_1$ is the $1$-form $\Theta^1$, also explicitly presented 
in \S 3. We pointed out that every constant multiple of $E_1$ is a Killing 
vector (equivalently, every constant multiple of $\Theta^1$ is a Killing 
covector) field of the Riemannian metric 
\begin{equation*} 
 \begin{aligned} 
  \tilde a :=& \ \epsilon^2 \, \Theta^1 \otimes \Theta^1 + 
                             \Theta^2 \otimes \Theta^2 + 
                             \Theta^3 \otimes \Theta^3        \\ 
            =& \ \omega^1 \otimes \omega^1 + 
                 \omega^2 \otimes \omega^2 + 
                 \omega^3 \otimes \omega^3 . 
 \end{aligned}  
\end{equation*} 
Motivated by this fact, let us stipulate the drift $1$-form to be 
$$\tilde b := \lambda \, \Theta^1 = \frac{\lambda}{\epsilon} \, \omega^1 , $$ 
where $\lambda$ is at the moment an arbitrary but nonzero constant. With 
respect to the metric $\tilde a$, this $\tilde b$ has constant Riemannian 
length 
$$\| \, \tilde b \, \| = \frac{ \, | \lambda | \, }{\epsilon} \ . $$

So the pair $\tilde a$ and $\tilde b$ do satisfy the first two criteria 
(among four) given by the Yasuda--Shimada theorem. Recall that those four  
criteria are necessary and sufficient for the Randers metric with data 
$\tilde a$, $\tilde b$ to have constant positive flag curvature $K$.  
The purpose of this section is to use the remaining two criteria to determine 
the constants $\epsilon$ and $\lambda$ in terms of $K$.  

We carry out our calculations in the moving frame $\{ e_1, e_2, e_3 \}$ and 
moving coframe $\{ \omega^1, \omega^2, \omega^3 \}$, which are 
$\tilde a$-orthonormal. The connection forms $\omega_q^{\ p}$ are displayed 
in \S 3. We have 
$$\tilde b_1 = \frac{\lambda}{\epsilon} \ , \ \ \ 
  \tilde b_2 = 0 , \ \ \ 
  \tilde b_3 = 0 . $$ 
And the relevant covariant differentiation formulas are: 
\begin{equation*} 
 \begin{aligned} 
  \tilde b_{p|q}   &= \left( d \, \tilde b_p  
                            - \tilde b_s \, \omega_p^{\ s} \right) (e_q) , \\ 
  \tilde b_{p|q|r} &= \left( d \, \tilde b_{p|q} 
                            - \tilde b_{s|q} \, \omega_p^{\ s}
                            - \tilde b_{p|s} \, \omega_q^{\ s} \right) (e_r) . 
 \end{aligned} 
\end{equation*}
Note that $\tilde b_{p|q}$ is skew-symmetric on its indices because 
$\tilde b$ is Killing. Consequently, the second covariant derivative 
$\tilde b_{p|q|r}$ is skew-symmetric in the indices $p$ and $q$. 
Straightforward calculations give:  
$$
\begin{pmatrix}
\, \tilde b_{1|1} & \tilde b_{1|2} & \tilde b_{1|3} \\ 
\, \tilde b_{2|1} & \tilde b_{2|2} & \tilde b_{2|3} \\ 
\, \tilde b_{3|1} & \tilde b_{3|2} & \tilde b_{3|3} 
\end{pmatrix} 
= 
\begin{pmatrix} 
0 & 0 & 0 \\ 
0 & 0 & - \lambda \\ 
0 & + \lambda & 0 
\end{pmatrix} \ ,  
$$
and 
\[
\begin{array}{lcllcllcl}   
 \tilde b_{1|2|1} & \= &                  0 , \ \  & 
 \tilde b_{1|3|1} & \= &                  0 , \ \  &  
 \tilde b_{2|3|1} & \= &                  0 ,              \\  
 \tilde b_{1|2|2} & \= & - \lambda \epsilon , \ \  & 
 \tilde b_{1|3|2} & \= &                  0 , \ \  &  
 \tilde b_{2|3|2} & \= &                  0 ,              \\  
 \tilde b_{1|2|3} & \= &                  0 , \ \  & 
 \tilde b_{1|3|3} & \= & - \lambda \epsilon , \ \  &  
 \tilde b_{2|3|3} & \= &                  0 .              
\end{array}
\]                          

Recall the third criterion given by Yasuda--Shimada. It reads 
$$\tilde b_{p|q|r} = T_{pqr} , \ \ \ \text{where} \ \ 
  T_{pqr} := K \, \left( \tilde a_{pr} \, \tilde b_q
                         - 
                         \tilde a_{qr} \, \tilde b_p \right) .$$
We find that  
\[
\begin{array}{lcllcllcl}   
 T_{121} & \= &                            0 , \ \  & 
 T_{131} & \= &                            0 , \ \  &  
 T_{231} & \= &                            0 ,              \\  
 T_{122} & \= & - \frac{\lambda K}{\epsilon} , \ \  & 
 T_{132} & \= &                            0 , \ \  &  
 T_{232} & \= &                            0 ,              \\  
 T_{123} & \= &                            0 , \ \  & 
 T_{133} & \= & - \frac{\lambda K}{\epsilon} , \ \  &  
 T_{233} & \= &                            0 .              
\end{array}
\]                          
Since the constant $\lambda$ is nonzero, the third Yasuda--Shimada criterion 
holds if and only if the dilation factor $\epsilon$ is given by 
$$\epsilon = \sqrt{K} . $$ 

We have just determined that, for our purpose, the Riemannian metric 
$\tilde a$ on $S^3$ should be 
\begin{equation*} 
 \begin{aligned} 
  \tilde a :=& \ K \, \Theta^1 \otimes \Theta^1 + 
                      \Theta^2 \otimes \Theta^2 + 
                      \Theta^3 \otimes \Theta^3        \\ 
            =& \ \omega^1 \otimes \omega^1 + 
                 \omega^2 \otimes \omega^2 + 
                 \omega^3 \otimes \omega^3 . 
 \end{aligned}  
\end{equation*} 
This information updates the Riemann curvature tensor we calculated (at the 
end of \S 3) to 
\[
\begin{array}{lcllcllcl}   
 \tilde R_{1212} & \= &    - K , \ \  & 
 \tilde R_{1213} & \= &      0 , \ \  &  
 \tilde R_{1223} & \= &      0 ,             \\  
 \tilde R_{1312} & \= &      0 , \ \  & 
 \tilde R_{1313} & \= &    - K , \ \  &  
 \tilde R_{1323} & \= &      0 ,             \\  
 \tilde R_{2312} & \= &      0 , \ \  & 
 \tilde R_{2313} & \= &      0 , \ \  &  
 \tilde R_{2323} & \= & 3K - 4 .              
\end{array}
\] 

Now examine the fourth (and last) criterion in the Yasuda--Shimada theorem. 
It requires the coupling between the drift $1$-form $\tilde b$ and the 
Riemannian metric $\tilde a$ to be such that  
$$\tilde R_{qprs} = \mathcal T_{qprs} ,$$ 
where $\mathcal T_{qprs}$ abbreviates the expression 
$$\begin{aligned}  
  &- K \left( 1 - \| \, \tilde b \, \|^2 \right) \, 
       \tilde a_{qr} \, \tilde a_{ps} 
   - K \left( \tilde a_{qr} \, \tilde b_p \, \tilde b_s 
              +  
              \tilde a_{ps} \, \tilde b_q \, \tilde b_r \right) 
   + \tilde b_{q|r} \, \tilde b_{p|s}                          \\ 
   &+ K \left( 1 - \| \, \tilde b \, \|^2 \right) \, 
        \tilde a_{qs} \, \tilde a_{pr} 
    + K \left( \tilde a_{qs} \, \tilde b_p \, \tilde b_r 
               +  
               \tilde a_{pr} \, \tilde b_q \, \tilde b_s \right) 
    - \tilde b_{q|s} \, \tilde b_{p|r}                          \\ 
   &+ 2 \, \tilde b_{q|p} \, \tilde b_{r|s} . 
       \end{aligned}$$ 
Calculations give 
\[
\begin{array}{lcllcllcl}   
 \mathcal T_{1212} & \= & - K , \ \  & 
 \mathcal T_{1213} & \= &   0 , \ \  &  
 \mathcal T_{1223} & \= &   0 ,             \\  
 \mathcal T_{1312} & \= &   0 , \ \  & 
 \mathcal T_{1313} & \= & - K , \ \  &  
 \mathcal T_{1323} & \= &   0 ,             \\  
 \mathcal T_{2312} & \= &   0 , \ \  & 
 \mathcal T_{2313} & \= &   0 , \ \  &  
 \mathcal T_{2323} & \= & 4 \lambda^2 - K .              
\end{array}
\] 
Therefore the fourth criterion of Yasuda--Shimada is satisfied if and only 
if 
$$\lambda = \pm \sqrt{ \, K-1 \, } =: \mathfrak s \sqrt{ \, K-1 \, } \ . $$
In particular, the constant positive flag curvature $K$ that one is striving 
for must be $\geqslant 1$. Anyway, the drift $1$-form has now been 
determined. It is 
$$\tilde b 
  := 
   \pm \sqrt{ \, K-1 \, } \, \Theta^1 
   = 
   \pm \sqrt{ \, \frac{K-1}{K} \, } \, \omega^1 . $$ 
It is somewhat amazing that the scaling multiple $\lambda$ (on $\tilde b$) 
asserts itself only at the $\tilde R_{2323} = \mathcal T_{2323}$ stage. 
 
When $K$ is equal to $1$, the scaling multiple $\lambda$ vanishes and the 
dilation factor $\epsilon$ is $1$. The Randers metric in question then reduces 
to the standard Riemannian one that $S^3$ inherits from Euclidean 
$\mathbb R^4$. This case is not of interest to us because we want 
non-Riemannian Randers spaces. So let us impose the restriction 
$$K > 1 . $$

\bigskip

\section{Two explicit descriptions of the resulting Randers metric} 

Let us recapitulate by giving explicit formulas for the Finsler functions 
of the Randers spaces we have just obtained. 

The first description is in terms of the non-holonomic frame 
$\{ e_1, e_2, e_3 \}$ and its coframe $\{ \omega^1, \omega^2, \omega^3 \}$, 
both globally defined on the manifold $M = S^3$. 
The frame is non-holonomic because each $e_p$ is a constant multiple of 
$E_p$, and the latter have nonzero Lie brackets amongst themselves (see \S 3). 

Anyway, the Riemannian metric is 
$$\tilde a = \omega^1 \otimes \omega^1 
             +
             \omega^2 \otimes \omega^2 
             +
             \omega^3 \otimes \omega^3 , $$  
and the drift $1$-form is 
$$\tilde b = \pm \sqrt{ \frac{ \, K-1 \, }{K} } \ \omega^1 . $$ 
Expanding arbitrary tangent vectors as 
$$y = \y^1 \, e_1 + \y^2 \, e_2 + \y^3 \, e_3 , $$ 
we see that the Finsler function is given by 
$$F(x, y) 
  = 
  \sqrt{ (\y^1)^2 + (\y^2)^2 + (\y^3)^2 } 
  \pm 
  \sqrt{ \frac{ \, K-1 \, }{K} } \ \y^1 . $$ 
It is remarkable that there is only an {\it implicit} dependence on the 
position $x$. We will revisit this deceptively simple formula in \S 7. 

The second description of our Randers metric is in terms of natural 
coordinates. For this purpose, we used the following parametrization of 
$S^3$: 
$$(x^o, x^1, x^2, x^3) 
  = 
  \frac{1}{\sqrt{ 1 + x^2 + y^2 + z^2 }} \ (c, x, y, z) , $$ 
where $c$ has the value $+1$ when dealing with the {\it right} 
hemisphere, and the value $-1$ when dealing with the {\it left} hemisphere. 
Obviously, the $y$ here no longer denotes generic tangent vectors on $S^3$. 
Instead, it is part of the collection $x$, $y$, $z$ that one typically uses 
when discussing Cartesian coordinates in Euclidean $\mathbb R^3$. 

This parametrization has the advantage that it imposes no restriction on 
$x$, $y$, $z$, and can be visualized as follows. Consider for example 
the right hemisphere of $S^3$, centered at the origin $O := (0, 0, 0, 0)$ 
of $\mathbb R^4$. Place a hyperplane $\mathbb R^3$ tangent to the sphere 
at the East Pole $(1, 0, 0, 0)$. Points on this tangent hyperplane are of 
the form $(1, x, y, z)$. Now we assign coordinates to any given position  
$P := (x^o, x^1, x^2, x^3)$ on the right hemisphere of $S^3$. Draw the 
straight line segment from the origin $O$ to $P$ and prolongate that until 
it intersects the tangent hyperplane. The point of intersection, being on 
that hyperplane, will have the form $(1, x, y, z)$. Just declare the natural 
coordinates of our point $P$ to be $(x, y, z)$. Note that points along the 
equator of the right hemisphere correspond to points at infinity on the 
tangent hyperplane. A similar story holds for the left hemisphere. 

The Riemannian metric $\tilde a$ and the drift $1$-form $\tilde b$ are 
$$\tilde a = \ K \, \Theta^1 \otimes \Theta^1 + 
                    \Theta^2 \otimes \Theta^2 + 
                    \Theta^3 \otimes \Theta^3  , $$
$$\tilde b = \pm \sqrt{ \, K-1 \, } \, \Theta^1 . $$
Simple computations give: 
\begin{equation*} 
 \begin{aligned} 
  \Theta^1 &= \frac{ \, + \, c dx - z dy + y dz \, } 
                   {1 + x^2 + y^2 + z^2} ,                \\ 
  \Theta^2 &= \frac{ \, + \, z dx + c dy - x dz \, } 
                   {1 + x^2 + y^2 + z^2} ,                \\ 
  \Theta^3 &= \frac{ \, - \, y dx + x dy + c dz \, } 
                   {1 + x^2 + y^2 + z^2} .            
 \end{aligned} 
\end{equation*} 
Denote generic tangent vectors on $S^3$ as 
$$u \frac{\partial}{\partial x} 
  + 
  v \frac{\partial}{\partial y} 
  +
  w \frac{\partial}{\partial z} \ . $$
Then the Finsler function for our Randers space is 
$$F(x,y,z;u,v,w) = \alpha(x,y,z;u,v,w) + \beta(x,y,z;u,v,w) , $$
with 
$$\alpha 
  = 
  \frac{ \, \sqrt{ K ( c u - z v + y w )^2 
                   + ( z u + c v - x w )^2 
                   + (-y u + x v + c w )^2 } \, } 
       { 1 + x^2 + y^2 + z^2 } \ , $$

$$\beta 
  = 
  \frac{ \, \pm \sqrt{ \, K-1 \, } \ ( c u - z v + y w ) \, }
       { \, 1 + x^2 + y^2 + z^2 \, } . $$      
As a reminder, $c = +1$ for the right hemisphere, and $c = -1$ 
for the left hemisphere. Some Maple programming will be carried out on  
this formula in the Appendix. 

The Riemannian norm of the drift $1$-form is 
$$\| \, \tilde b \, \| 
  = 
  \sqrt{ \frac{K-1}{K} } 
  = 
  \sqrt{ 1 - \frac{1}{K} } \ . $$
It is given by this value whether one is using the $\tilde a$-orthonormal 
frame or natural coordinates. Since $K \geqslant 1$ (note: $K = 1$ 
corresponds to the standard sphere $S^3$), we see that 
$\| \, \tilde b \, \|$ is strictly less than $1$. This is both necessary and 
sufficient (\cite{AIM}, \cite{BCS}) for $F$ to be positive and strongly 
convex on the entire slit tangent bundle $TM \smallsetminus 0$.  

\bigskip

\section{Spray curvatures and flag curvatures} 

The flag curvature is defined much like the sectional curvature of 
Riemannian geometry. One begins with a connection, differentiates 
appropriately, and then performs some standard contractions. Let us first give 
a synopsis of the procedure. 

The connection in use (from a family that contains the handful of name brand 
ones) gives rise to a $hh$-curvature tensor $R^{\ i}_{j \ kl}$, 
which in the Riemannian case is precisely the familiar Riemann curvature. 
The precise formula of $R^{\ i}_{j \ kl}$ in terms of, say, the Chern 
connection, has been given a pedagogical treatment in \cite{BCS}. Since that 
specific formula does not concern us in this paper, we shall omit it. 

At any point $x$ on $M$, one associates with each flagpole---a nonzero 
vector $y$ in $T_xM$---the {\it spray} curvature   
$$K^i_{\ k} := y^j \, R^{\ i}_{j \ kl} \ y^l . $$
For those familiar with \cite{BCS}, this $K^i_{\ k}$ is $F^2$ times something 
called the `predecessor of the flag curvature.' It is actually the first 
among two spray curvatures (\cite{B3}, \cite{D}). The second one, using the 
version described in \cite{D} and multiplying by $F$, is equal to the 
$hv$-curvature ${}^b P^{\ i}_{j \ kl}$ of the Berwald connection, but is 
irrelevant to the purpose of this paper. 

Our notation $K^i_{\ k}$ is the same as Berwald's \cite{B1}. Rund \cite{R} 
uses $H^i_{\ k}$ instead. The spray curvature is robust enough that it 
does not depend on our choice of connection. However, the above definition is 
often not an efficient way to compute the spray curvature because one 
has to obtain the tensor $R^{\ i}_{j \ kl}$ first. 

\begin{quote} 
A word of caution about notation. For us, $K_{ik}$ shall simply mean 
$g_{ij} K^j_{\ k}$. It can be shown (see \cite{BCS}) that $K_{ik}$ is 
symmetric in $i$, $k$, and $y^i K_{ik}$, $K_{ik} y^k$ both vanish. Rund did  
{\it not} subscribe to this system. He defined his $H_{ik}$ by first tracing 
the indices $j$ and $l$ on $\third [(H^j_{\ l})_{y^k} - (H^j_{\ k})_{y^l}]$, 
forming an intermediate $H_k$, and then taking the $y$-partial $(H_k)_{y^i}$.  
See pages 129 and then 125 of \cite{R}. As a result, his $H_{ik}$ was not 
symmetric and $y^i H_{ik}$ did not vanish for him. 
\end{quote}

The flagpole $y$, together with any vector $V \in T_xM$ transversal to $y$, 
specifies a flag based at the point $x \in M$. The {\it flag} curvature of 
the resulting flag is the quantity 
$$K(y,V)
  := 
  \frac{V^i \, K_{ik} \, V^k}
       { \, g(y,y) g(V,V) - [g(y,V)]^2 \, } \ ,$$
where $g := g_{ij} dx^i \otimes dx^j$ is the fundamental tensor defined 
in \S 2. The Finsler metric is said to have constant flag curvature if 
$K(y,V)$ has the same constant value $K$ for any choice of $y$ and $V$. 

Specialize now to the case of constant flag curvature $K$. 
Euler's theorem for homogeneous functions implies that $g(y,y) = F^2(x,y)$. 
Let's substitute this into the above formula and rearrange it into the form 
$$V^i \, K_{ik} \, V^k 
  = 
  K \, F^2 \, \{ g(V,V) - [g(\ell,V)]^2 \} , $$
where $\ell$ means $y/F$. As such, this statement makes formal sense even if 
$V$ is not transversal to $y$. Since both sides are given by symmetric 
bilinear forms, a standard polarization identity gives 
$$U^i \, K_{ik} \, W^k 
  = 
  K \, F^2 \, \{ g(U,W) - g(\ell,U) g(\ell,W) \} . $$
In particular, 
$$K_{ik} = K \, F^2 \, \left( g_{ik} - F_{y^i} F_{y^k} \right) . $$  
Here, we have used the identity $g_{sj} \ell^s = F_{y^j}$, again a 
consequence of Euler's theorem. Raising the index $i$ with the inverse 
matrix $g^{ij}$ of the fundamental tensor, we get 
$$K^i_{\ k} 
  = 
  K \, F^2 \, \left( \delta^i_{\ k} - \frac{y^i}{F} F_{y^k} \right) . $$
We shall use this criterion as the characterization of constant flag 
curvature. 

Next, we describe a computationally friendlier way to access the spray 
curvature $K^i_{\ k}$. The history behind this better approach goes all the 
way back to Berwald \cite{B1} in his study of path spaces and projective 
geometry. See Rund \cite{R} for a careful exposition, and references therein. 

One begins with the geodesic spray coefficients 
$$G^i := \half \, \gamma^i_{\ jk} \, y^j y^k , $$
where 
$$\gamma^i_{\ jk} 
  :=  
  g^{il} \ \frac{1}{2} 
             \left( \frac{\partial g_{lj}}{\partial x^k} 
                    -  
                    \frac{\partial g_{jk}}{\partial x^l} 
                    + 
                    \frac{\partial g_{kl}}{\partial x^j} \right) $$ 
are the fundamental tensor's formal Christoffel symbols of the second 
kind. Our definition of $G^i$ here is $\half$ times that given in the book 
\cite{BCS}. As a result, the present definition agrees exactly with the 
one used in \cite{AIM} and \cite{R}. Chapter IV formula (6.3) in \cite{R} 
tells us that 
$$K^i_{\ k} 
  = 
  2 (G^i)_{x^k} 
  - 
  y^j (G^i)_{x^j y^k} 
  - 
  (G^i)_{y^j} (G^j)_{y^k} 
  + 
  2 \, G^j (G^i)_{y^j y^k} . $$
For ease of exposition, let us refer to this as Berwald's formula. It is an 
instructive exercise to manipulate the original expression for the 
spray curvature, namely $y^j R_{j \ kl}^{\ i} \, y^l$, into the one 
we just displayed. See \S 12.5B in \cite{BCS} if guidance is needed.  

The above formula expresses the spray curvature in terms of partial 
derivatives of the geodesic spray coefficients. As such, it is quite useful 
for machine computations of the spray curvature. The Maple codes that 
implement this sort of computation are given in the Appendix. 

The said formula is also invaluable for calculating spray curvatures by 
hand. However, its efficient implementation for this purpose often requires 
some additional maneuvers as described in \cite{S1}. To get us set up for 
\S 7, we specialize this technique to the Randers setting and present it here. 

The geodesic spray coefficients of every Randers metric is equal to that 
of the underlying Riemannian metric, plus a perturbation term. Symbolically, 
we write 
$$G^i = \tilde G^i + \zeta^i . $$
Here, 
$$\tilde G^i := \half \, \tilde \gamma^i_{\ jk} \, y^j y^k , $$
where 
$$\tilde \gamma^i_{\ jk}
  := 
  \tilde a^{il} \ \frac{1}{2} 
             \left( \frac{\partial \tilde a_{lj}}{\partial x^k} 
                    -  
                    \frac{\partial \tilde a_{jk}}{\partial x^l} 
                    + 
                    \frac{\partial \tilde a_{kl}}{\partial x^j} \right) $$  
are the Christoffel symbols of the second kind of the Riemannian metric 
$\tilde a$. The specific form of the perturbation term $\zeta^i$ has been 
worked out in \cite{BCS}, but shall not concern us until \S 7. 

Substitute the above $G^i$ into Berwald's formula for the spray 
curvature. After covariantizing (to be explained below) all $x$-partial 
derivatives of $\zeta^i$ and using Euler's theorem for homogeneous functions 
whenever appropriate, we get 
$$K^i_{\ k} 
  = 
  \tilde K^i_{\ k} 
  + 
  \left \{ 2 \zeta^i_{\ |k} 
           - 
           y^j (\zeta^i_{\ |j})_{y^k} 
           - 
           (\zeta^i)_{y^j} (\zeta^j)_{y^k} 
           + 
           2 \, \zeta^j (\zeta^i)_{y^j y^k} 
  \right \} . $$ 
This formula will be our centerpiece in \S 7. 
The first term on the right-hand side, $\tilde K^i_{\ k}$, is the spray 
curvature of the Riemannian metric $\tilde a$. It can readily be  
calculated in the $S^3$ examples because the Riemann curvature 
$\tilde R_{jikl}$ of $\tilde a$ is already at our disposal (\S 4). 

The remaining group of terms on the right-hand side involves horizontal 
covariant differentiation on the slit tangent bundle $TM \smallsetminus 0$. 
Let us explain how those are carried out. The geodesic spray coefficients 
$\tilde G^i$ are first used to define 
$$\tilde N^i_{\ j} := (\tilde G^i)_{y^j}$$
and then 
$$\frac{\delta}{\delta x^j} 
  := 
  \frac{\partial}{\partial x^j} 
  - 
  \tilde N^l_{\ j} \, \frac{\partial}{\partial y^l} \ .$$
The vector fields $\frac{\delta}{\delta x^j}$ are declared to be  
horizontal---the horizontal lift of $\frac{\partial}{\partial x^j}$ 
to $TM \smallsetminus 0$. And the $\tilde N^i_{\ j}$ are said to have 
produced a nonlinear Ehresmann connection on the slit tangent bundle. In 
order to achieve the sought horizontal covariant differentiation on 
$\zeta^i$, we horizontally lift the Levi-Civita (Christoffel) connection of 
$\tilde a$, and then let the resulting object act on $\zeta^i$. This may seem 
abstract but is operationally quite simple:  
$$\zeta^i_{\ |j} 
  := 
  \frac{\delta}{\delta x^j} \zeta^i 
  \, + \, 
  \zeta^l \, \tilde \gamma^i_{\ lj} \ .$$

Incidentally, note that 
$$\tilde \gamma^i_{\ jk} = ( \tilde G^i )_{y^j y^k} .$$
This complements nicely with the fact that the nonlinear connection 
is given by the first $y$-partial derivatives of $\tilde G^i$. Thus 
$$(\zeta^i)_{x^j} 
  = 
  \zeta^i_{\ |j} 
  + 
  (\zeta^i)_{y^l} (\tilde G^l)_{y^j} 
  - 
  \zeta^l (\tilde G^i)_{y^l y^j} .$$
The use of this formula is what we mean above by covariantizing the 
$x$-partial derivatives of $\zeta^i$. We did that to make each term inside 
the expression $\{ \cdots \}$ manifestly tensorial. Of course, the quantity 
$\{ \cdots \}$ {\it as a whole} is already tensorial, whether we 
covariantize or not. There is no need to covariantize the $y$-partial 
derivatives because they already transform tensorially (see \cite{BCS} for an 
exposition). 

\bigskip

\section{Direct proof of constant positive flag curvature} 

As we mentioned in the Introduction, the two published proofs (\cite{YS} and 
\cite{M}) of the Yasuda--Shimada theorem both contain arguments that we 
find difficult to follow. Given that, it is prudent to prove directly that 
our Randers metrics have constant positive flag curvature $K$. The numerical 
evidence documented in the Appendix helps sustain such a bare-hands' proof. 
Strategically, we compute both sides of the constant curvature criterion 
$$K^p_{\ r} 
  = 
  K \, F^2 \, 
  \left( \delta^p_{\ r} - \frac{\y^p}{F} F_{\y^r} \right)$$  
and ascertain that they are indeed equal. The calculations rely on Berwald's 
formula, and will be done with the help of an $\tilde a$-orthonormal frame 
instead of natural coordinates. The resulting proof is geometrical and is 
distinctly different from the approaches used in \cite{YS}, \cite{M}. It 
presents a new perspective on the subject, and is worthy of further 
development.

\subsection {The right-hand side of the constant curvature criterion} 

\ \smallskip  

Let us begin with the formula for $F$ in terms of the $\tilde a$-orthonormal 
frame $\{ e_1, e_2, e_3 \}$ on $S^3$. This was first given in \S 5. 
For ease of presentation, we relabel the tangent space coordinates 
$\y^1$, $\y^2$, $\y^3$ as $\u$, $\v$, $\w$, respectively. 
These coordinates arise from the expansion of arbitrary tangent vectors $y$ 
in terms of the basis $\{ e_1, e_2, e_3 \}$. The cosmetically altered formula 
for $F$ now reads 
$$F(x, y) 
  = 
  \sqrt{ \u^2 + \v^2 + \w^2 } 
  +  
  \mathfrak s \sqrt{ \frac{ \, K-1 \, }{K} } \ \u 
  =: 
  \alpha + \beta \ .$$ 
Here, $\mathfrak s = \pm 1$ keeps track of the two choices of the drift 
$1$-form. See the end of \S 4. These choices of sign have nothing to do with 
the right or left hemisphere of $S^3$. The latter are associated with the 
parameter $c$ of \S 5. 

Almost by inspection, we have 
\begin{equation*} 
 \begin{aligned} 
  F_\u &= \frac{\u}{\alpha} + \mathfrak s \sqrt{ \frac{K-1}{K} } \ ,     \\ 
  F_\v &= \frac{\v}{\alpha} \ ,                                          \\ 
  F_\w &= \frac{\w}{\alpha} \ . 
 \end{aligned} 
\end{equation*} 
Using these, the right-hand side 
$$K \, \tau^p_{\ r} 
  := 
  K \, F^2 \, 
  \left( \delta^p_{\ r} - \frac{\y^p}{F} F_{\y^r} \right)$$
of our constant curvature criterion is readily computed. We find that
$$
\begin{pmatrix} 
K \, \tau^{1}_{\ 1} & K \, \tau^{1}_{\ 2} & K \, \tau^{1}_{\ 3} \\ 
K \, \tau^{2}_{\ 1} & K \, \tau^{2}_{\ 2} & K \, \tau^{2}_{\ 3} \\ 
K \, \tau^{3}_{\ 1} & K \, \tau^{3}_{\ 2} & K \, \tau^{3}_{\ 3}
\end{pmatrix} 
$$ 
is equal to 
$$
\begin{pmatrix}
  \frac{F}{\alpha} \, K (\v^2 + \w^2) & 
- \frac{F}{\alpha} \, K \u\v          & 
- \frac{F}{\alpha} \, K \u\w                                               \\ 
- \frac{F}{\alpha} \, K \u\v  
  - F \, K \v \, \mathfrak s \sqrt{ \frac{K-1}{K} } & 
  F \, K (F - \frac{\v^2}{\alpha})                  & 
- \frac{F}{\alpha} \, K \v\w                                               \\ 
- \frac{F}{\alpha} \, K \u\w  
  - F \, K \w \, \mathfrak s \sqrt{ \frac{K-1}{K} } & 
- \frac{F}{\alpha} \, K \v\w                        &
  F \, K (F - \frac{\w^2}{\alpha})        
\end{pmatrix} \ . 
$$

\subsection {Perturbation terms in the geodesic spray coefficients} 

\ \smallskip 

The left-hand side of our constant curvature criterion is the spray 
curvature $K^p_{\ r}$. These nine components can be calculated using 
Berwald's formula. To this end, we start with the geodesic spray coefficients 
$G^i$. As we mentioned in \S 6, the ones used in this paper, {\it unlike} 
those of \cite{BCS}, already have the factor of $\half$ built in. It is shown 
in the same reference that for Randers metrics, 
$$G^i = \tilde G^i + \zeta^i ,$$
where $\tilde G^i := \half \, \tilde \gamma^i_{\ jk} \, y^j y^k$ are the 
geodesic spray coefficients of the underlying Riemannian metric $\tilde a$. 

The perturbation terms $\zeta^i$ transform like the components of a tensor. 
They do so because 
$$\zeta^i 
  = 
  \half \left( \gamma^i_{\ jk} - \tilde \gamma^i_{\ jk} \right) y^j y^k ,$$
showing that $\zeta$ arises from the difference of two connections. More 
abstractly, consider the vector bundle $TM$ that sits over $M$. Using the 
projection map $\pi : TM \smallsetminus 0 \rightarrow M$, we can pull that 
back to obtain a vector bundle $\pi^*TM$ that sits over the manifold 
$TM \smallsetminus 0$. This just means that over each point 
$(x, y) \in TM \smallsetminus 0$, we have erected a copy of $T_xM$. Our 
$\zeta^i$ then transforms like a section of this pulled-back vector bundle. 

The coordinate bases $\{ \partial_{x^i} \}$, the Riemannian metric 
$\tilde a$, together with the field of $\tilde a$-orthonormal frames  
$\{ e_p \}$, can all be transplanted to the fibres of $\pi^*TM$ 
(note: {\it not} to $TM \smallsetminus 0$). As on $M$, the transplants of 
$\{ \partial_{x^i} \}$ and $\{ e_p \}$ are related through the $3$-bein 
$u_p^{\ i}$ and its matrix inverse $v^p_{\ i}$: 
$$e_p = u_p^{\ i} \, \partial_{x^i} \ , 
  \ \ \ \ 
  \omega^p = v^p_{\ i} \, dx^i . $$
Each $u_p^{\ i}$ and $v^p_{\ i}$ is a function on $M$. In other words, 
they depend on $x$ only. 

The tensorial property mentioned above allows us to transform the 
coordinate description $\zeta^i$ into the orthonormal frame description 
$$\zeta^p := v^p_{\ i} \, \zeta^i .$$
With this in mind, and the factor of $\half$ emphasized earlier, an 
expression derived in \cite{BCS} is transcribed. It gives   
$$\zeta^p 
  := 
  \half \, \alpha \, \tilde b_{q|r}   
  \left \{ 
           (\tilde a^{pq} \y^r - \tilde a^{pr} \y^q)  
           + 
           \ell^p \, ( \y^q \, \tilde b^r 
                       - 
                       \y^r \, \tilde b^q )   
           +  
           \frac{\ell^p \, \y^q \y^r}{\alpha}   
  \right \} .$$ 
Here, $\ell^p := \y^p / F$ is to be distinguished from 
$\tilde \ell^p := \y^p / \alpha$. We shall show that these $\zeta^p$ 
components are actually quite simple. 

In \S 4, we determined the intermediate constants $\epsilon$ and $\lambda$. 
These can be used to update the Levi-Civita (Christoffel) connection 
$\omega_q^{\ p}$ of the Riemannian metric $\tilde a$, the components of the 
$1$-form $\tilde b$, and its covariant derivatives $\tilde b_{q|r}$. 
A quick glance at \S 3 and \S 4 tells us that the updated versions read: 
$$
\begin{pmatrix}
\omega_1^{\ 1} & \omega_1^{\ 2} & \omega_1^{\ 3} \\ 
\omega_2^{\ 1} & \omega_2^{\ 2} & \omega_2^{\ 3} \\ 
\omega_3^{\ 1} & \omega_3^{\ 2} & \omega_3^{\ 3} 
\end{pmatrix} 
= 
\begin{pmatrix} 
0  & - \sqrt{K} \omega^3 & \sqrt{K} \omega^2 \\ 
\sqrt{K} \omega^3 & 0 & \left(\sqrt{K}-\frac{2}{\sqrt{K}}\right) \omega^1 \\ 
- \sqrt{K} \omega^2 & \left(\frac{2}{\sqrt{K}}-\sqrt{K}\right) \omega^1 & 0 
\end{pmatrix} \ , 
$$
together with  
$$\tilde b_1 = \mathfrak s \sqrt{ \frac{ \, K-1 \, }{K} } \ , \ \ \ 
  \tilde b_2 = 0 , \ \ \ 
  \tilde b_3 = 0 ,$$ 
and 
$$
\begin{pmatrix}
\, \tilde b_{1|1} & \tilde b_{1|2} & \tilde b_{1|3} \\ 
\, \tilde b_{2|1} & \tilde b_{2|2} & \tilde b_{2|3} \\ 
\, \tilde b_{3|1} & \tilde b_{3|2} & \tilde b_{3|3} 
\end{pmatrix} 
= 
\begin{pmatrix} 
0 & 0                        & 0                        \\ 
0 & 0                        & - \mathfrak s \sqrt{K-1} \\ 
0 & + \mathfrak s \sqrt{K-1} & 0 
\end{pmatrix} \ .  
$$

Since the present $\tilde b_{q|r}$ is skew-symmetric and most of its 
components are zero, the perturbation term $\zeta^p$ simplifies drastically 
to: 
$$\zeta^p 
  = 
  \alpha \, \mathfrak s \sqrt{K-1} 
  \left \{ \tilde a^{p3} \y^2 - \tilde a^{p2} \y^3 
           + 
           \ell^p (\y^3 \, \tilde b^2 - \y^2 \tilde b^3) 
  \right \} . $$
Furthermore, we are in an $\tilde a$-orthonormal frame, so 
$\tilde b^r := \tilde a^{rs} \tilde b_s = \delta^{rs} \tilde b_s$ 
is numerically the same as $\tilde b_r$. Given that $\tilde b_2$ and 
$\tilde b_3$ are both zero, the terms multiplying $\ell^p$ drop out, and 
$$\zeta^p 
  = 
  \alpha \, \mathfrak s \sqrt{K-1} 
  \left \{ \delta^{p3} \y^2 - \delta^{p2} \y^3 \right \} . $$ 
In other words, 
\begin{equation*}
 \begin{aligned} 
  \zeta^1 &= 0 ,                                                       \\ 
  \zeta^2 &= - \, \mathfrak s \sqrt{ \, K-1 \, } \ \w \, \alpha ,      \\ 
  \zeta^3 &= + \, \mathfrak s \sqrt{ \, K-1 \, } \ \v \, \alpha . 
 \end{aligned} 
\end{equation*} 
Recall that in the $\tilde a$-orthonormal frame description, the quantity 
$\alpha$ is simply $\sqrt{\u^2 + \v^2 + \w^2}$. So the components $\zeta^p$ 
contain no explicit dependence on $x$. However, $\u$, $\v$, $\w$ are $\y^p$, 
which can be expressed as $v^p_{\ i} \, y^i$. Thus the dependence on $x$ is 
there, albeit implicitly through $v^p_{\ i}$.

\subsection {The spray curvature of the Riemannian metric} 

\ \smallskip 

In \S 6, we described a refined understanding of Berwald's formula. It is 
a case study of the general strategy discussed in \cite{S1}. The essence 
is that given the decomposition $G^i = \tilde G^i + \zeta^i$, the spray 
curvature splits as well: 
$$K^i_{\ k} 
  = 
  \tilde K^i_{\ k} 
  + 
  \left \{ 2 \zeta^i_{\ |k} 
           - 
           y^j (\zeta^i_{\ |j})_{y^k} 
           - 
           (\zeta^i)_{y^j} (\zeta^j)_{y^k} 
           + 
           2 \, \zeta^j (\zeta^i)_{y^j y^k} 
  \right \} . $$ 
Each side of this equation transforms like a section of the pulled-back 
tensor bundle $\pi^*TM \otimes \pi^*T^*M$, which sits over 
$TM \smallsetminus 0$. 

Contract both sides of the above equation with $v^p_{\ i} \, u_r^{\ k}$.  
This converts its expression from the coordinate to the 
$\tilde a$-orthonormal perspective: 
$$K^p_{\ r} 
  = 
  \tilde K^p_{\ r} 
  + 
  \left \{ 2 \zeta^p_{\ |r} 
           - 
           \y^q (\zeta^p_{\ |q})_{\y^r} 
           - 
           (\zeta^p)_{\y^q} (\zeta^q)_{\y^r} 
           + 
           2 \, \zeta^q (\zeta^p)_{\y^q \y^r} 
  \right \} . $$ 
The first term on the right-hand side is the spray curvature of the 
Riemannian metric $\tilde a$. Since the Riemann curvature tensor 
$\tilde R_{qprs}$ of $\tilde a$ has been determined in \S 4, we can use that 
information to compute the said spray curvature: 
$$\tilde K^p_{\ r} = \y^q \, \tilde R^{\ p}_{q \ rs} \ \y^s . $$ 
Keep in mind that we are in an $\tilde a$-orthonormal frame, so the numerical 
values of tensor components are unaffected by raising or lowering any  
index, as long as it is done using $\tilde a$. Fairly routine calculations 
tell us that the nine spray curvatures 
$$
\begin{pmatrix} 
\, \tilde K^1_{\ 1} & \tilde K^1_{\ 2} & \tilde K^1_{\ 3} \, \\ 
\, \tilde K^2_{\ 1} & \tilde K^2_{\ 2} & \tilde K^2_{\ 3} \, \\ 
\, \tilde K^3_{\ 1} & \tilde K^3_{\ 2} & \tilde K^3_{\ 3} \, 
\end{pmatrix} 
$$
are equal to 
$$
\begin{pmatrix} 
K (\v^2 + \w^2) & - K \u\v               & - K \u\w                \\ 
- K \u\v        & K \u^2 + (4 - 3K) \w^2 & (3K - 4) \v\w           \\
- K \u\w        & (3K - 4) \v\w          & K \u^2 + (4 - 3K) \v^2
\end{pmatrix} . 
$$
The symmetry here is expected because 
$\tilde K^p_{\ r} := \tilde a^{ps} \tilde K_{sr}$ has the same 
numerical value as $\tilde K_{pr}$, and the latter is symmetric by general 
principles. In contradistinction to this, the indices on $K^p_{\ r}$ are 
raised and lowered with $g$ rather than $\tilde a$. Since the basis we 
are using on the fibres of $\pi^*TM$ is not $g$-orthonormal, we do {\it not} 
expect $K^p_{\ r}$ to be numerically equal to the symmetric $K_{pr}$.

\subsection {A simplified formula for horizontal covariant derivatives} 

\ \smallskip 

In the refined version of Berwald's formula, we have the terms 
$$2 \zeta^p_{\ |r} 
  - 
  \y^q (\zeta^p_{\ |q})_{\y^r} 
  - 
  (\zeta^p)_{\y^q} (\zeta^q)_{\y^r} 
  + 
  2 \, \zeta^q (\zeta^p)_{\y^q \y^r} . $$ 
These involve both the $\y$-partial derivatives of $\zeta^p$ as well 
as its horizontal covariant derivatives. The purpose of this subsection 
is to derive a simplified formula for these horizontal covariant derivatives, 
one that will facilitate their computation in \S 7.5. 

In \S 6, we explained how the horizontal covariant derivatives are to be 
carried out in the context of natural coordinates. Namely, 
$$\zeta^i_{\ |k} 
  = 
  \frac{\delta}{\delta x^k} \zeta^i 
  \, + \, 
  \zeta^l \, \tilde \gamma^i_{\ lk} \ .$$
Note that $\tilde \gamma^i_{\ lj} \, dx^j$ are the connection forms with 
respect to the coordinate basis. When these connection forms are pulled back 
to $TM \smallsetminus 0$ and evaluated on the horizontal lifts 
(of $\partial_{x^k}$ to $TM \smallsetminus 0$) 
$$\frac{\delta}{\delta x^k} 
  := 
  \frac{\partial}{\partial x^k} 
  - 
  \tilde N^i_{\ k} \, \frac{\partial}{\partial y^i} \ ,$$
they see only the $\partial_{x^k}$ part but not the $\partial_{y^i}$ 
part. The above expression for $\zeta^i_{\ |k}$ can then be recast as 
$$\zeta^i_{\ |k} 
  = 
  \left \{ d \zeta^i 
           \, + \, 
           \zeta^l \, (\tilde \gamma^i_{\ lj} \, dx^j) \right \}  
  \left( \frac{\delta}{\delta x^k} \right) \ .$$ 

In the current subsection, we are using the basis $\{ e_p \otimes \omega^r \}$ 
(instead of the coordinate one) on the pulled-back bundle 
$\pi^*TM \otimes \pi^*T^*M$. The analogous formula for $\zeta^p_{\ |r}$ is 
$$\zeta^p_{\ |r} 
  = 
  \left \{ d \zeta^p 
           \, + \, 
           \zeta^s \, \omega_s^{\ p} \right \}  
  ( \hat e_r ) \ ,$$
where $\omega_s^{\ p}$ are the connection forms displayed in \S 7.2. The  
vector field $\hat e_r$ is the horizontal lift of $e_r$, and has the formula 
$$\hat e_r 
  := 
  e_r 
  - 
  \y^s \, \omega_s^{\ p} \left( e_r \right) \, \partial_{\y^p} .$$       
This is not unreasonable, given that in the horizontal lift of 
$\partial_{x^k}$, the nonlinear connection can be re-expressed as follows: 
$$\tilde N^i_{\  k} 
  = 
  \tilde \gamma^i_{\ lk} \, y^l 
  = 
  y^l \, \left( \tilde \gamma^i_{\ lj} \, dx^j \right) 
  \left( \partial_{x^k} \right) . $$

Note, however, that we have yet to define $e_r$, and for that matter 
$\partial_{x^k}$, as objects on $TM \smallsetminus 0$. Here are the 
definitions: 
\begin{itemize} 
\item The $1$-forms $dx^k$ are pulled-back to $TM \smallsetminus 0$ and given 
      the same name. They are then complemented by the $dy^k$ (with $y^k$ 
      coming from $y^k \partial_{x^k}$) to form a coordinate basis for the 
      cotangent bundle of $TM \smallsetminus 0$. The natural dual of this 
      basis is denoted $\{ \partial_{x^k}, \partial_{y^k} \}$. The objects in 
      this basis are local vector fields on $TM \smallsetminus 0$. They are 
      the ones that enter the definition of $\frac{\delta}{\delta x^k}$. It 
      is simply an abuse of notation to employ the same symbol 
      $\partial_{x^k}$ for the vector fields here and those on $M$. 
\smallskip 
\item In the same spirit, we pull back the $1$-forms $\omega^r$ to 
      $TM \smallsetminus 0$ and retain the same name. These are complemented 
      by the $d \y^r$ (with $\y^r$ arising from $\y^r e_r$) to give a 
      non-holonomic basis. Denote the natural dual of this basis as 
      $\{ e_r, \partial_{\y^r} \}$. Just like the $\partial_{x^k}$ case, we 
      are abusing the notation when we call these $e_r$ and those on $M$ by 
      the same name. Nonetheless, we now know what $e_r$ means on 
      $TM \smallsetminus 0$, and the formula for the horizontal lift 
      $\hat e_r$ makes good sense.  
\end{itemize} 
Do not confuse the $e_r$ and $\partial_{x^k}$ we just defined with the 
transplants (see \S 7.2) that live on the pulled-back vector bundle 
$\pi^*TM$, even though they share the same notation merely for the sake of 
economy. 

One can't help but wonder how $e_r$ is related to $\partial_{x^k}$ as vector 
fields on $TM \smallsetminus 0$. To find out, take an arbitrary differentiable 
function $f$ on $TM \smallsetminus 0$ and let $e_r$ act on it. By the chain 
rule, 
$$(df)(e_r) 
  = 
  \left \{ (\partial_{x^k}f) \, dx^k 
           + 
           (\partial_{y^k}f) \, dy^k \right \} (e_r) .$$
\begin{itemize} 
\item Note that $y^k = \y^s u_s^{\ k}$ gives 
      $dy^k = u_s^{\ k} d \y^s + \y^s d u_s^{\ k}$. 
      It is to minimize confusion in this statement that we have chosen 
      the font $\y^s$ throughout the paper, when dealing with 
      orthonormal expansions of $y$. 
\item The relation $dx^k = u_q^{\ k} \omega^q$ is valid on 
      $TM \smallsetminus 0$ because it holds on $M$ and the pull-back 
      operation preserves algebraic statements. 
\item Since $\{ e_p, \partial_{\y^p} \}$ is by definition the  
      natural dual of $\{ \omega^p, d \y^p \}$, 
      one must have $\omega^q (e_r) = \delta^q_{\ r}$ and $d \y^s (e_r) = 0$. 
\end{itemize}
Substituting these observations into our chain rule statement, and removing 
the test function $f$ afterwards, we obtain the somewhat surprising formula 
$$e_r 
  = 
  u_r^{\ k} \, \partial_{x^k} 
  + 
  \y^s \, du_s^{\ k}(e_r) \, \partial_{y^k} 
  \ \ \ 
  \text{on} \ TM \smallsetminus 0 . $$
In contradistinction to that, a similar calculation gives 
$$\partial_{\y^r} = u_r^{\ k} \, \partial_{y^k} 
  \ \ \ 
  \text{on} \ TM \smallsetminus 0 . $$
Here, one does have to invoke $(d u_s^{\ k})(\partial_{\y^r}) = 0$. This 
holds because $u_s^{\ k}$ is a function on $M$, so its differential is a 
linear combination of $dx^j$ and hence of $\omega^q$. Such a relation is 
preserved under pull-back, and $\omega^q (\partial_{\y^r}) = 0$. 

Let us return to the horizontal covariant derivative 
$$\zeta^p_{\ |r} 
  = 
  \left \{ d \zeta^p 
           \, + \, 
           \zeta^s \, \omega_s^{\ p} \right \} (\hat e_r) \ ,$$
with $\hat e_r := e_r - 
      \y^s \, \omega_s^{\ t} (e_r) \, \partial_{\y^t}$.        
This can stand two more reductions before we put it to use in \S 7.5. 
\begin{itemize} 
\item On $M$, the connection forms $\omega_s^{\ p}$ are linear combinations 
      of the $\omega^q$, which are in turn linear combinations of the $dx^i$. 
      This remains so under the pull-back to $TM \smallsetminus 0$, 
      so $\omega_s^{\ p} (\partial_{y^k}) = 0$. 
      Since $\partial_{\y^t} = u_t^{\ k} \, \partial_{y^k}$, we must have 
      $\omega_s^{\ p} (\partial_{\y^t}) = 0$ as well.
      Therefore 
      $$\omega_s^{\ p} (\hat e_r) = \omega_s^{\ p} (e_r) . $$ 
\item At the end of \S 7.2, we deduced that the components $\zeta^p$ depend  
      only on $\u$, $\v$, $\w$, namely $\y^q$. Formally taking the 
      differential yields $d \zeta^p = (\zeta^p)_{\y^q} d \y^q$. On the other 
      hand, linear algebra says that 
      $d \zeta^p =   (d \zeta^p)(e_q) \, \omega^q 
                   + (d \zeta^p)(\partial_{\y^q}) \, d \y^q$. Comparing the 
      two statements gives $(d \zeta^p)(e_q) = 0$, which in turn allows us 
      to conclude that 
      $$(d \zeta^p)(\hat e_r) 
        = 
        - \y^s \, \omega_s^{\ t} (e_r) \, (\zeta^p)_{\y^t} .$$
\end{itemize} 
The two reductions above lead to the following: 
$$\zeta^p_{\ |r} 
  = 
  \left \{ \zeta^s \, \omega_s^{\ p} 
           - 
           (\zeta^p)_{\y^t} \ \y^s \, \omega_s^{\ t} \right \} (e_r) . $$
We shall use this formula in \S 7.5 to compute horizontal covariant 
derivatives. By the way, in the context of \cite{BCS}, the $\omega_q^{\ p}$ 
here happens to be the Chern connection forms of the Riemannian metric 
$\tilde a$.

\subsection {Derivatives of the perturbation terms} 

\ \smallskip 

In this subsection, we tabulate the requisite $\y^r$-partial derivatives 
of $\zeta^p$. Then we use the concluding formula of \S 7.4, together with 
$\omega_q^{\ p}$ from \S 7.2 and $\omega^s (e_r) = \delta^s_{\ r}$, to 
calculate the horizontal covariant derivatives that we need.  Recall that 
\begin{equation*}
 \zeta^1 = 0 , \ \ \                                                        
 \zeta^2 = - \, \mathfrak s \sqrt{ \, K-1 \, } \ \w \, \alpha , \ \ \  
 \zeta^3 = + \, \mathfrak s \sqrt{ \, K-1 \, } \ \v \, \alpha ,  
\end{equation*} 
with 
$$\alpha := \sqrt{ \, \u^2 + \v^2 + \w^2 \, } . $$
Here, $\u = \y^1$, $\v = \y^2$, $\w = \y^3$ and, from \S 7.1, 
$\mathfrak s = \pm 1$. Thus all $\y^r$-partials of $\zeta^1$ are zero. 
On the other hand, $\zeta^1_{\ |r}$ do not necessarily vanish because 
$\zeta^2$ and $\zeta^3$ are involved in the calculations. 

Let us list the results. We introduce some abbreviations for ubiquitous 
quantities in order to reduce clutter. 

$$\text{With} \ \ 
  \diamondsuit := \frac{\mathfrak s}{ \, \alpha \, } \ \sqrt{ \, K-1 \, } : $$
\[
\begin{array}{lcllcl}   
 (\zeta^2)_{\y^1} & \= & 
        - \diamondsuit \, (\u \w)            \ \  & 
 (\zeta^3)_{\y^1} & \= &                
        + \diamondsuit \, (\u \v)                  \\ 
 (\zeta^2)_{\y^2} & \= &                
        - \diamondsuit \, (\v \w)            \ \  &  
 (\zeta^3)_{\y^2} & \= &     
        + \diamondsuit \, (\u^2+2\v^2+\w^2)        \\ 
 (\zeta^2)_{\y^3} & \= &                
        - \diamondsuit \, (\u^2+\v^2+2\w^2)  \ \  &    
 (\zeta^3)_{\y^3} & \= &                
        + \diamondsuit \, (\v \w) 
\end{array}
\] 

\bigskip

$$\text{With} \ \ 
  \bigstar := \frac{\mathfrak s}{ \, \alpha^3 \, } \ \sqrt{ \, K-1 \, } : $$
\[
\begin{array}{lcllcl}   
 (\zeta^2)_{\y^1 \y^2} & \= & 
        + \bigstar \, (\u \v \w)              & 
 (\zeta^3)_{\y^1 \y^2} & \= &                
        + \bigstar \, \u (\u^2+\w^2)                  \\ 
 (\zeta^2)_{\y^1 \y^3} & \= &                
        - \bigstar \, \u (\u^2+\v^2)          &  
 (\zeta^3)_{\y^1 \y^3} & \= &     
        - \bigstar \, (\u \v \w)                      \\ 
 (\zeta^2)_{\y^2 \y^2} & \= &                
        - \bigstar \, \w (\u^2+\w^2)          &  
 (\zeta^3)_{\y^2 \y^2} & \= &     
        + \bigstar \, \v (3\u^2+2\v^2+3\w^2)          \\ 
 (\zeta^2)_{\y^2 \y^3} & \= &                
        - \bigstar \, \v (\u^2+\v^2)          &  
 (\zeta^3)_{\y^2 \y^3} & \= &     
        + \bigstar \, \w (\u^2+\w^2)                  \\ 
 (\zeta^2)_{\y^3 \y^3} & \= &                
        - \bigstar \, \w (3\u^2+3\v^2+2\w^2)  &    
 (\zeta^3)_{\y^3 \y^3} & \= &                
        + \bigstar \, \v (\u^2+\v^2) 
\end{array}
\] 

\bigskip

$$\text{With} \ \ 
  \heartsuit 
  := 
  \alpha \ \mathfrak s \ \sqrt{ \, K-1 \, } \ \sqrt{ \, K \, } : $$
\[
\begin{array}{lcllcllcl}   
 \zeta^1_{\ |1} & \= & 0                  \ \ \ \ \ \ & 
 \zeta^2_{\ |1} & \= & 0                    \ \ \ \ \ &
 \zeta^3_{\ |1} & \= & 0                               \\ 
 \zeta^1_{\ |2} & \= & - \heartsuit \, \v \ \ \ \ \ \ & 
 \zeta^2_{\ |2} & \= & + \heartsuit \, \u   \ \ \ \ \ &
 \zeta^3_{\ |2} & \= & 0                               \\ 
 \zeta^1_{\ |3} & \= & - \heartsuit \, \w \ \ \ \ \ \ & 
 \zeta^2_{\ |3} & \= & 0                    \ \ \ \ \ &
 \zeta^3_{\ |3} & \= & + \heartsuit \, \u              
\end{array}
\] 

\bigskip 

\noindent Lastly, 

$$\text{With} \ \ 
  \clubsuit  
  := 
  \frac{1}{\alpha} \ \mathfrak s \ \sqrt{ \, K-1 \, } \ \sqrt{ \, K \, } : $$
\[ 
\begin{array}{lcllcl} 
 \left(\zeta^2_{\ |2}\right)_{\y^1} & \= & 
        + \clubsuit \, (2\u^2+\v^2+\w^2)           \ \  & 
 \left(\zeta^3_{\ |3}\right)_{\y^1} & \= &                
        + \clubsuit \, (2\u^2+\v^2+\w^2)                  \\ 
 \left(\zeta^2_{\ |2}\right)_{\y^2} & \= &                
        + \clubsuit \, (\u \v)                     \ \  &  
 \left(\zeta^3_{\ |3}\right)_{\y^2} & \= &     
        + \clubsuit \, (\u \v)                            \\ 
 \left(\zeta^2_{\ |2}\right)_{\y^3} & \= &                
        + \clubsuit \, (\u \w)                     \ \  &    
 \left(\zeta^3_{\ |3}\right)_{\y^3} & \= &                
        + \clubsuit \, (\u \w) 
\end{array}
\]  
These two columns are equal because $(\zeta^2_{\ |2}) = (\zeta^3_{\ |3})$.  
Also, we have 
\[ 
\begin{array}{lcllcl} 
 \left(\zeta^1_{\ |2}\right)_{\y^1} & \= & 
        - \clubsuit \, (\u \v)                     \ \  & 
 \left(\zeta^1_{\ |3}\right)_{\y^1} & \= &                
        - \clubsuit \, (\u \w)                            \\ 
 \left(\zeta^1_{\ |2}\right)_{\y^2} & \= &                
        - \clubsuit \, (\u^2+2\v^2+\w^2)           \ \  &  
 \left(\zeta^1_{\ |3}\right)_{\y^2} & \= &     
        - \clubsuit \, (\v \w)                            \\ 
 \left(\zeta^1_{\ |2}\right)_{\y^3} & \= &                
        - \clubsuit \, (\v \w)                     \ \  &    
 \left(\zeta^1_{\ |3}\right)_{\y^3} & \= &                
        - \clubsuit \, (\u^2+\v^2+2\w^2) 
\end{array}
\]

\subsection {The spray curvature of our Randers metric} 

\ \smallskip 

According to \S 7.3, the spray curvature of the Randers metric 
$$F := \alpha + \beta$$ 
has the structure 
$$K^p_{\ r} = \tilde K^p_{\ r} + \mathcal E^p_{\ r} \ ,$$ 
where 
$$\mathcal E^p_{\ r}
  := 
  2 \zeta^p_{\ |r} 
  - 
  \y^q (\zeta^p_{\ |q})_{\y^r} 
  - 
  (\zeta^p)_{\y^q} (\zeta^q)_{\y^r} 
  + 
  2 \, \zeta^q (\zeta^p)_{\y^q \y^r} . $$ 
The term $\tilde K^p_{\ r}$ represents the nine spray curvatures of the 
underlying Riemannian metric $\tilde a$. We have already calculated those 
in \S 7.3. 

Using the derivatives listed in \S 7.5, the above quantities 
$\mathcal E^p_{\ r}$ are computed. We find that: 

\begin{equation*}
 \begin{aligned} 
  \mathcal E^1_{\ 1}   
  &= 
  + \clubsuit \, \u (\v^2+\w^2) ,                                 \\
  \mathcal E^1_{\ 2} 
  &= 
  - \clubsuit \, \v (\u^2) ,                                      \\ 
  \mathcal E^1_{\ 3} 
  &= 
  - \clubsuit \, \w (\u^2) ;                                      \\
  \mathcal E^2_{\ 1} 
  &= 
  - \clubsuit \, \v (2\u^2+\v^2+\w^2)  - (K-1) \u \v ,            \\ 
  \mathcal E^2_{\ 2} 
  &= 
  + \clubsuit \, \u (2\u^2+\v^2+2\w^2) + (K-1) (\u^2+4\w^2) ,     \\ 
  \mathcal E^2_{\ 3} 
  &= 
  - \clubsuit \, \u \v \w - 4(K-1) \v \w ;                        \\
  \mathcal E^3_{\ 1} 
  &= 
  - \clubsuit \, \w (2\u^2+\v^2+\w^2) - (K-1) \u \w ,             \\ 
  \mathcal E^3_{\ 2} 
  &= 
  - \clubsuit \, \u \v \w  - 4(K-1) \v \w ,                       \\ 
  \mathcal E^3_{\ 3} 
  &= 
  + \clubsuit \, \u (2\u^2+2\v^2+\w^2) + (K-1) (\u^2+4\v^2) .     
 \end{aligned} 
\end{equation*} 

Adding these $\mathcal E^p_{\ r}$ to the corresponding $\tilde K^p_{\ r}$ 
of \S 7.3, we obtain the spray curvatures $K^p_{\ r}$. After that, it 
is a matter of routine algebra to make sure each $K^p_{\ r}$ is equal to the 
$K \tau^p_{\ r}$ calculated in \S 7.1. Such is indeed the case. This means 
that the Randers metric in question satisfies the criterion 
$$K^p_{\ r} 
  = 
  K \, F^2 \, 
  \left( \delta^p_{\ r} - \frac{\y^p}{F} F_{\y^r} \right) . $$
Therefore it has constant flag curvature $K$ and our proof is complete.   

\bigskip

\section{Discussion} 

Given any Finsler manifold $(M, F)$, its geodesics are paths in $M$, just 
like the Riemannian case. The defining equation for geodesics with constant 
Finslerian speed is 
$$\ddot x^i + \gamma^i_{\ jk} \, \dot x^j \dot x^k = 0 . $$
In terms of the coefficients $G^i$ of \S 6, this system 
of second order quasi-linear equations can be rewritten as 
$$\ddot x^i + 2 \, G^i = 0 . $$
This is one reason the $G^i$ are called the geodesic spray coefficients. 

A Finsler manifold is said to be projectively flat if $M$ can be covered 
by privileged coordinate charts in which the above quasi-linear system 
of ordinary differential equations becomes a linear system. 
Projective flatness has been characterized by Douglas \cite{D} in his study 
of path spaces. When specialized to the Finsler setting, his result states 
that: 
\begin{quote} 
A Finsler manifold $(M, F)$ of dimension $\geqslant 3$ is projectively flat 
if and only if its projective Weyl and Douglas tensors both vanish.
\end{quote} 
A slightly modified characterization, due to Berwald \cite{B2}, applies to 
dimension $2$. There, the vanishing of the Weyl tensor is replaced by 
another criterion. See page 144 of Rund \cite{R}, or \cite{S1}, for a 
detailed account and further references. 

The projective Weyl tensor is to be distinguished from the conformal Weyl 
tensor in Riemannian geometry. The latter is not defined in dimension $2$ and 
always vanishes in dimension $3$. It turns out that the vanishing of the 
projective Weyl tensor (which has four indices) is equivalent to the 
vanishing of its reduced version $W^i_{\ k}$. This is explained on pages 141 
and 142 of \cite{R}. 

The tensor $W^i_{\ k}$ is defined in terms of the spray curvature and its 
derivatives. To that end, one first constructs the scalar quantity 
$$\mathfrak K := \frac{1}{n-1} \, K^i_{\ i} \, .$$
Comparing with \S 7.6 of \cite{BCS}, we see that this scalar $\mathfrak K$ 
can be interpreted as $F^2$ times the average of $n-1$ appropriately chosen 
flag curvatures. Using $\mathfrak K$, the reduced projective Weyl tensor is 
defined as 
$$W^i_{\ k} 
  := 
  K^i_{\ k} 
  - 
  \mathfrak K \, \delta^i_{\ k} 
  - 
  \frac{1}{n+1} \ y^i 
  \left \{ (K^j_{\ k})_{y^j} - (\mathfrak K)_{y^k} \right \} . $$

A result of Matsumoto's (see \cite{AIM}) says that the vanishing of 
$W^i_{\ k}$ is both necessary and sufficient for the flag curvature to 
depend at most on the position $x$ and the flagpole $y \in T_xM$, but not on 
the transverse edges (denoted $V$ in \S 6). In particular, since our Randers 
metrics on $S^3$ have constant flag curvature $K$, we deduce that $W^i_{\ k}$ 
must vanish for all these examples. This can also be verified by a direct 
computation. 

The Douglas tensor $D^i_{\ jkl}$ (denoted $B^i_{\ jkl}$ in \cite{R}) is 
defined as follows: 
\begin{equation*} 
 \begin{aligned} 
  D^i_{\ jkl} 
  := 
  \ \ \ \ 
  &(G^i)_{y^j y^k y^l}                                       \\ 
   - \ 
  &\frac{1}{n+1} 
   \left \{   \delta^i_{\ j} (G^h)_{y^h y^k y^l} 
            + \delta^i_{\ k} (G^h)_{y^h y^l y^j} 
            + \delta^i_{\ l} (G^h)_{y^h y^j y^k} \right \}    \\ 
   - \ 
  &\frac{1}{n+1} \ y^i \, (G^h)_{y^h y^j y^k y^l} . 
 \end{aligned} 
\end{equation*} 
In \cite{BM}, B\'asc\'o and Matsumoto proved that a Randers space has 
vanishing Douglas tensor if and only if the drift $1$-form $\tilde b$ is 
closed. For our examples, $\tilde b$ is a constant multiple of $\Theta^1$, 
see \S 4. Since $d \Theta^1 = 2 \Theta^2 \wedge \Theta^3 \not= 0$, we 
conclude that the Douglas tensor does not vanish for our Randers metrics. 

Now we know that each member among our family of Randers metrics on $S^3$ 
has vanishing projective Weyl tensor but nonzero Douglas tensor. Therefore,  
by Douglas' theorem, none of them can be projectively flat.

\bigskip

\appendix

\section{Some numerical evidence} 

The purpose of this Appendix is to apply Maple to the Randers metric 
$$F(x,y,z;u,v,w) = \alpha(x,y,z;u,v,w) + \beta(x,y,z;u,v,w)$$
in natural coordinates, with 
$$\alpha 
  = 
  \frac{ \, \sqrt{ K ( c u - z v + y w )^2 
                   + ( z u + c v - x w )^2 
                   + (-y u + x v + c w )^2 } \, } 
       { 1 + x^2 + y^2 + z^2 } \ , $$

$$\beta 
  = 
  \frac{ \, \pm \sqrt{ \, K-1 \, } \ ( c u - z v + y w ) \, }
       { \, 1 + x^2 + y^2 + z^2 \, } . $$      
Recall that $c = +1$ for the right hemisphere of $S^3$, and $c = -1$ 
for the left hemisphere. 

In Finsler geometry, it is not uncommon for simple formulae to quickly 
mushroom into unmanageable expressions. We find that machine computations 
have consistently extracted useful information and insights, to the point 
that meaningful followup questions can be asked. We value every opportunity 
to cultivate this synergy between Finsler geometry and modern computing. 
It is hoped that by producing the Maple codes here, we are initiating other 
geometers into a fruitful aspect of experimental mathematics.  

Here is our plan: 
\begin{itemize} 
\item We first use Maple to calculate the spray curvature {\it \`a la} 
      Berwald: 
      $$K^i_{\ k} 
        = 
        2 (G^i)_{x^k} 
        - 
        y^j (G^i)_{x^j y^k} 
        - 
        (G^i)_{y^j} (G^j)_{y^k} 
        + 
        2 \, G^j (G^i)_{y^j y^k} . $$
\item After that, we ask Maple to check whether $F$ satisfies the 
      characterization of having constant flag curvature $K$: 
      $$K^i_{\ k} 
        = 
        K \, F^2 \, \left( \delta^i_{\ k} - \frac{y^i}{F} F_{y^k} \right) . $$ 
\end{itemize}

\subsection{The Finsler function in natural coordinates} 

\ \smallskip 

$>$ P:=c*u+y*w--z*v; 

$>$ Q:=c*v+z*u--x*w; 

$>$ R:=c*w+x*v--y*u; 
 
Let us restrict our attention to the right hemisphere of $S^3$. So 

$>$ c:=+1; 

\noindent Define, for lack of a better name: 

$>$ den:=1+x\^2+y\^2+z\^2;

\noindent Then  

$>$ alpha:=(1/den)*sqrt(K*P\^2+Q\^2+R\^2); 

$>$ beta:=(1/den)*sqrt(K--1)*P; 
 
\noindent Here, we have simply chosen the $+$ sign in the drift term $\beta$. 
A moment's thought shows that there is no loss of generality in doing so. 
The Finsler function $F$ and its associated Lagrangian $L$ are: 

$>$ F:=alpha+beta; 
 
$>$ L:=(1/2)*F\^2; 

\noindent Use of the semi-colon instructs Maple to display the input formulas.

\subsection{The covariant form of the geodesic spray coefficients} 

\ \smallskip 

By that, we are referring to the quantities
$$G_i 
  :=  
  \half \, \gamma_{ijk} \, y^j y^k 
  = 
  \fourth  
  \left( g_{ij, x^k} - g_{jk, x^i} + g_{ki, x^j} \right) y^j y^k . $$
Raising the index with the inverse of the fundamental tensor $g_{ij}$ gives 
the contravariant form $G^i$. That will be carried out in \S A.4. In terms of 
the Lagrangian $L$, it is not difficult to show that 
$$G_i 
  = 
  \half \left( L_{y^i x^j} \, y^j - L_{x^i} \right) . $$

\noindent Our natural coordinates $x^1$, $x^2$, $x^3$ are denoted $x$, $y$, 
$z$, and the induced tangent space coordinates $y^1$, $y^2$, $y^3$ are 
denoted $u$, $v$, $w$. The covariant form of the geodesic spray coefficients 
are named (by us) $g1$, $g2$, $g3$ in Maple. Thus 

$>$ Lx:=diff(L,x): 

$>$ Ly:=diff(L,y): 
  
$>$ Lz:=diff(L,z): 
 
$>$ g1:=(1/2)*(u*diff(Lx,u)+v*diff(Ly,u)+w*diff(Lz,u)--Lx):
  
$>$ g2:=(1/2)*(u*diff(Lx,v)+v*diff(Ly,v)+w*diff(Lz,v)--Ly):
  
$>$ g3:=(1/2)*(u*diff(Lx,w)+v*diff(Ly,w)+w*diff(Lz,w)--Lz):
 
\noindent Use of the colon instructs Maple to suppress the computed answers.

\subsection{The inverse of the fundamental tensor} 

\ \smallskip 

The inverse $g^{ij}$ of the fundamental tensor $g_{ij}$ (\S 2) is needed in 
order to obtain the contravariant form $G^i$ of the geodesic spray 
coefficients. This inverse has a 
standard formula  
$$g^{ij} 
  = 
  \frac{\alpha}{F} \, \tilde a^{ij} 
  +  
  \frac{\alpha^2}{F^2} \,  
  \frac{ \, \beta + \alpha \, \| \tilde b \|^2 \, }{F} \,  
  \tilde \ell^i \tilde \ell^j                    
  -  
  \frac{\alpha^2}{F^2} \, \left( \tilde \ell^i \, \tilde b^j 
                                 +  
                                 \tilde \ell^j \, \tilde b^i \right) , $$
where $\tilde \ell^i := \frac{y^i}{\alpha}$, 
$\tilde b^i :=  \tilde a^{ij} \tilde b_j$, and $\tilde a^{ij}$ is the 
inverse of the Riemannian metric $\tilde a_{ij}$. For a pedagogical 
derivation, one can consult \cite{BCS}. 

In our natural coordinates, the matrix of the Riemannian metric $\tilde a$ 
is
$$
\frac{1}{ \text{den}^2 }
\begin{pmatrix}
K + z^2 + y^2   & - Kcz + cz - xy & Kcy - xz - cy   \\ 
- Kcz + cz - xy & Kz^2 + 1 + x^2  & - Kyz           \\ 
Kcy - xz - cy   & - Kyz           & Ky^2 + x^2 + 1 
\end{pmatrix} ,$$
with $den:=1+x^2+y^2+z^2$. The inverse of $\tilde a_{ij}$ has been computed 
elsewhere using Maple. We assign it the name $Aij$ in our codes:  
\begin{itemize}
\item[$>$] A11:=(1/K)*((x\^2+1)*(x\^2+K*z\^2+K*y\^2+1)):
\item[$>$] A12:=(1/K)*(x\^3*y+x*y\^3*K+x*y--c*z*x\^2--c*z
                 \newline ${\ }$ \ \ \ \ \ \ \  
                 +K*c*z*x\^2+K*c*z+K*y*z\^2*x):
\item[$>$] A13:=(1/K)*(K*y\^2*z*x+x\^3*z+x*z\^3*K+x*z
                \newline ${\ }$ \ \ \ \ \ \ \  
                --K*c*y*x\^2--K*c*y+c*y*x\^2+c*y):
\item[$>$] A21:=A12:
\item[$>$] A22:=(1/K)*(y\^2*x\^2+K*x\^2+(K--1)*2*x*y*c*z
                \newline ${\ }$ \ \ \ \ \ \ \  
                +y\^4*K+z\^2+K+y\^2*K*z\^2+2*K*y\^2):
\item[$>$] A23:=(1/K)*(x\^2*z*y+(K--1)*c*z\^2*x--(K--1)*x*y\^2*c
                \newline ${\ }$ \ \ \ \ \ \ \  
                +K*y\^3*z+K*y*z\^3--z*y+2*y*K*z):
\item[$>$] A31:=A13:
\item[$>$] A32:=A23:
\item[$>$] A33:=(1/K)*(z\^2*x\^2+K*x\^2--(K--1)*2*x*y*c*z
                \newline ${\ }$ \ \ \ \ \ \ \  
                +y\^2*K*z\^2+z\^4*K+2*K*z\^2+K+y\^2):
\end{itemize}  
As we have previously decided, $c = +1$ because we want to focus on the right 
hemisphere. The (covariant) components $\tilde b_i$ of our drift $1$-form are 
$$\frac{\sqrt{K-1}}{ \, 1 + x^2 + y^2 + z^2 \, } \, 
  \left( \, c, \, -z, \, +y \, \right) . $$
Denote its contravariant components $\tilde b^i := \tilde a^{ij} \tilde b_j$ 
in Maple as $B1$, $B2$, $B3$. On the right hemisphere: 

$>$ kappa:=sqrt(K--1)/den:
 
$>$ B1:=kappa*(A11--z*A12+y*A13):
  
$>$ B2:=kappa*(A21--z*A22+y*A23):
  
$>$ B3:=kappa*(A31--z*A32+y*A33):
 
\noindent In Maple, let us denote $\tilde \ell^i := \frac{y^i}{\alpha}$ as 
$tel1$, $tel2$, $tel3$, and the Riemannian norm 
$$\| \, \tilde b \, \| 
  = 
  \sqrt{ \frac{K-1}{K} } 
  = 
  \sqrt{ 1 - \frac{1}{K} } \ \ \text{as} \ B .$$
We have:  

$>$ tel1:=u/alpha:
  
$>$ tel2:=v/alpha:
  
$>$ tel3:=w/alpha:
  
$>$ B:=sqrt(1--(1/K)):
 
\noindent To reduce clutter, let us also introduce two abbreviations: 

$>$ rho:=alpha/F:

$>$ phi:=(beta+alpha*B\^2)/F:
 
\noindent Now we are ready to give the formula for the inverse $g^{ij}$ of the 
fundamental tensor. Denote that inverse, in Maple, as $Gij$. Then: 
 
$>$ G11:=rho*A11+rho\^2*phi*tel1*tel1--rho\^2*(tel1*B1+tel1*B1):
  
$>$ G12:=rho*A12+rho\^2*phi*tel1*tel2--rho\^2*(tel1*B2+tel2*B1):
  
$>$ G13:=rho*A13+rho\^2*phi*tel1*tel3--rho\^2*(tel1*B3+tel3*B1):
  
$>$ G21:=G12:
  
$>$ G22:=rho*A22+rho\^2*phi*tel2*tel2--rho\^2*(tel2*B2+tel2*B2):
  
$>$ G23:=rho*A23+rho\^2*phi*tel2*tel3--rho\^2*(tel2*B3+tel3*B2):
  
$>$ G31:=G13:
  
$>$ G32:=G23:
  
$>$ G33:=rho*A33+rho\^2*phi*tel3*tel3--rho\^2*(tel3*B3+tel3*B3):

\subsection{The contravariant form of the geodesic spray coefficients} 

\ \smallskip 

These are obtained by taking the covariant form $G_i$ of the coefficients 
and raising the index with the inverse $g^{ij}$ of the fundamental tensor. 

$>$ G1:=G11*g1+G12*g2+G13*g3:
  
$>$ G2:=G21*g1+G22*g2+G23*g3:
  
$>$ G3:=G31*g1+G32*g2+G33*g3:

\subsection{Getting set up for Berwald's formula} 

\ \smallskip 

We assign names to the first and second order partial derivatives of $G^i$. 
This will avoid having to re-compute them every time they are needed. 
Our Maple codes for Berwald's formula should run more efficiently as a result 
of this move. 

$>$ G1x:=diff(G1,x):

$>$ G1ux:=diff(G1u,x): 

$>$ G1vx:=diff(G1v,x):

$>$ G1wx:=diff(G1w,x):

\noindent Duplicate the above with $x$ replaced by $y$, $z$, $u$, $v$, $w$. 
 
\noindent Then duplicate all with $G1$ replaced successively 
          by $G2$ and $G3$.

\subsection{The spray curvature {\it \`a la} Berwald's formula} 

\ \smallskip 

We now get Maple to calculate the nine components $K^i_{\ k}$ of the 
spray curvature using Berwald's formula 
$$K^i_{\ k} 
  = 
  2 (G^i)_{x^k} 
  - 
  y^j (G^i)_{x^j y^k} 
  - 
  (G^i)_{y^j} (G^j)_{y^k} 
  + 
  2 \, G^j (G^i)_{y^j y^k} . $$
Even though the covariant form $K_{ik}$ of the spray curvature is 
symmetric in its two indices, the same {\it cannot} usually be said of the 
type $\binom11$ form $K^i_{\ k}:=g^{ij} K_{jk}$. So, even at a purely 
numerical level, $K^k_{\ i} \not= K^i_{\ k}$ in general, unless of course one 
is using a $g$-orthonormal frame. 

Let us use $\mathit{Kay}[i,k]$ as the Maple names for those nine spray 
curvatures. The reason for not using $K[i,k]$ is because of the next command. 

$>$ Kay:=array(1..3,1..3): 

\noindent Had we used $K[i,k]$, the above would read 
``K:=array(1..3,1..3):'', which would wreak havoc with our codes 
since the Maple variable $K$ already stands for something else (namely the 
constant positive flag curvature). 

\begin{itemize}
\item[$>$] Kay[1,1]:=2*G1x
                     \newline ${\ }$ \ \ \ \ \ \ \ \ \ \ \ \  
                     --G1u*G1u--G1v*G2u--G1w*G3u
                     \newline ${\ }$ \ \ \ \ \ \ \ \ \ \ \ \  
                     --u*G1ux--v*G1uy--w*G1uz
                     \newline ${\ }$ \ \ \ \ \ \ \ \ \ \ \ \  
                     +2*G1*G1uu+2*G2*G1uv+2*G3*G1uw:
\item[$>$] Kay[1,2]:=2*G1y
                     \newline ${\ }$ \ \ \ \ \ \ \ \ \ \ \ \  
                     --G1u*G1v--G1v*G2v--G1w*G3v
                     \newline ${\ }$ \ \ \ \ \ \ \ \ \ \ \ \  
                     --u*G1vx--v*G1vy--w*G1vz
                     \newline ${\ }$ \ \ \ \ \ \ \ \ \ \ \ \  
                     +2*G1*G1vu+2*G2*G1vv+2*G3*G1vw:
\item[$>$] Kay[1,3]:=2*G1z
                     \newline ${\ }$ \ \ \ \ \ \ \ \ \ \ \ \  
                     --G1u*G1w--G1v*G2w--G1w*G3w
                     \newline ${\ }$ \ \ \ \ \ \ \ \ \ \ \ \  
                     --u*G1wx--v*G1wy--w*G1wz
                     \newline ${\ }$ \ \ \ \ \ \ \ \ \ \ \ \  
                     +2*G1*G1wu+2*G2*G1wv+2*G3*G1ww:
\item[$>$] Kay[2,1]:=2*G2x
                     \newline ${\ }$ \ \ \ \ \ \ \ \ \ \ \ \  
                     --G2u*G1u--G2v*G2u--G2w*G3u
                     \newline ${\ }$ \ \ \ \ \ \ \ \ \ \ \ \  
                     --u*G2ux--v*G2uy--w*G2uz
                     \newline ${\ }$ \ \ \ \ \ \ \ \ \ \ \ \  
                     +2*G1*G2uu+2*G2*G2uv+2*G3*G2uw:
\item[$>$] Kay[2,2]:=2*G2y
                     \newline ${\ }$ \ \ \ \ \ \ \ \ \ \ \ \  
                     --G2u*G1v--G2v*G2v--G2w*G3v
                     \newline ${\ }$ \ \ \ \ \ \ \ \ \ \ \ \  
                     --u*G2vx--v*G2vy--w*G2vz
                     \newline ${\ }$ \ \ \ \ \ \ \ \ \ \ \ \  
                     +2*G1*G2vu+2*G2*G2vv+2*G3*G2vw:
\item[$>$] Kay[2,3]:=2*G2z
                     \newline ${\ }$ \ \ \ \ \ \ \ \ \ \ \ \  
                     --G2u*G1w--G2v*G2w--G2w*G3w
                     \newline ${\ }$ \ \ \ \ \ \ \ \ \ \ \ \  
                     --u*G2wx--v*G2wy--w*G2wz
                     \newline ${\ }$ \ \ \ \ \ \ \ \ \ \ \ \ 
                     +2*G1*G2wu+2*G2*G2wv+2*G3*G2ww:
\item[$>$] Kay[3,1]:=2*G3x 
                     \newline ${\ }$ \ \ \ \ \ \ \ \ \ \ \ \  
                     --G3u*G1u--G3v*G2u--G3w*G3u
                     \newline ${\ }$ \ \ \ \ \ \ \ \ \ \ \ \  
                     --u*G3ux--v*G3uy--w*G3uz
                     \newline ${\ }$ \ \ \ \ \ \ \ \ \ \ \ \  
                     +2*G1*G3uu+2*G2*G3uv+2*G3*G3uw:
\item[$>$] Kay[3,2]:=2*G3y 
                     \newline ${\ }$ \ \ \ \ \ \ \ \ \ \ \ \   
                     --G3u*G1v--G3v*G2v--G3w*G3v
                     \newline ${\ }$ \ \ \ \ \ \ \ \ \ \ \ \  
                     --u*G3vx--v*G3vy--w*G3vz
                     \newline ${\ }$ \ \ \ \ \ \ \ \ \ \ \ \  
                     +2*G1*G3vu+2*G2*G3vv+2*G3*G3vw:
\item[$>$] Kay[3,3]:=2*G3z
                     \newline ${\ }$ \ \ \ \ \ \ \ \ \ \ \ \  
                     --G3u*G1w--G3v*G2w--G3w*G3w
                     \newline ${\ }$ \ \ \ \ \ \ \ \ \ \ \ \   
                     --u*G3wx--v*G3wy--w*G3wz
                     \newline ${\ }$ \ \ \ \ \ \ \ \ \ \ \ \  
                     +2*G1*G3wu+2*G2*G3wv+2*G3*G3ww:
\end{itemize}

\subsection{The criterion for having constant flag curvature}

\ \smallskip 

Finally, we ask Maple to check whether our Randers metric $F$ has constant 
flag curvature $K$. The criterion we shall use has been derived in \S 6.  
It reads  
$$K^i_{\ k} 
  = 
  K \, F^2 \, \left( \delta^i_{\ k} - \frac{y^i}{F} F_{y^k} \right) 
  =: K \, \tau^i_{\ k} , $$
where $\tau^i_{\ k}$ equals $F^2$ times the terms inside the parentheses. The 
pertinent Maple codes are: 

$>$ Fu:=diff(F,u):
  
$>$ Fv:=diff(F,v):
  
$>$ Fw:=diff(F,w):

$>$ tau:=array(1..3,1..3):
 
$>$ tau[1,1]:=F\^2--u*F*Fu:
  
$>$ tau[1,2]:=    --u*F*Fv:
  
$>$ tau[1,3]:=    --u*F*Fw:
  
$>$ tau[2,1]:=    --v*F*Fu:
  
$>$ tau[2,2]:=F\^2--v*F*Fv:
  
$>$ tau[2,3]:=    --v*F*Fw:
  
$>$ tau[3,1]:=    --w*F*Fu:
  
$>$ tau[3,2]:=    --w*F*Fv:
  
$>$ tau[3,3]:=F\^2--w*F*Fw:
 
\noindent The condition we are striving for is $K^i_{\ k} = K \tau^i_{\ k}$. 
To see if that's the case, we form the difference of the two sides, 
and also their ratio for good measure. 

$>$ with(linalg):
 
$>$ f:= (i,k) --$>$ Kay[i,k]--K*tau[i,k]:
 
$>$ dif:=matrix(3,3,f): 
   
\noindent Hopefully the answers are zero.   
 
$>$ h:= (i,k) --$>$ Kay[i,k]/(K*tau[i,k]): 
 
$>$ quot:=matrix(3,3,h):
     
\noindent These quotients should all be 1.

\subsection{The verdict} 

\ \smallskip 

It was no problem for Maple to calculate {\it symbolically} those differences 
and quotients defined in \S A.7. However, attempts to get Maple to 
symbolically simplify the answers to a $0$ or a $1$ consistently crashed 
because the resulting expressions were too large. 

Given that, we did the next best thing. We randomly selected numerical values 
for the position coordinates $x$, $y$, $z$, the velocity variables $u$, $v$, 
$w$, and the positive flag curvature $K$. Then we asked Maple to evaluate 
those differences and quotients at the stipulated $x$, $y$, $z$, $u$, $v$, 
$w$, and $K$. In retrospect, it now appears that for the purpose of comparing 
$K^i_{\ k}$ to $K \tau^i_{\ k}$, forming the ratio of the two terms works 
better than taking their difference. 

We noticed that the amount of RAM used in our computations is routinely 
in excess of 1 Gigabyte. It is not clear whether this is attributable to 
any inefficiency in our codes. A sampling of our numerical results are 
as follows: 
 
$>$ simplify(eval(dif[1,1],[x=1.0,y=2.0,z=--3.0,
                       
\qquad\qquad\qquad\qquad\qquad\ \ u=3.1416,v=2.78,w=137.0,K=29.0])); 
                              
\qquad\qquad                                 .009000000000
 
\smallskip

$>$ factor(simplify(eval(quot[1,1],[x=1.0,y=2.0,z=--3.0,
                      
\qquad\qquad\qquad\qquad\qquad\qquad\qquad\ u=3.1416,v=2.78,w=137.0,K=29.0])));
 
\qquad\qquad                                  1.000000007

\smallskip 
 
$>$ simplify(eval(dif[1,2],[x=9.0,y=7.0,z=--5.0,

\qquad\qquad\qquad\qquad\qquad\ \ u=3.1416,v=2.78,w=137.0,K=31.0]));
 
\qquad\qquad                                   .008218000000
 
\smallskip

$>$ factor(simplify(eval(quot[1,2],[x=9.0,y=7.0,z=--5.0,
                              
\qquad\qquad\qquad\qquad\qquad\qquad\qquad\ u=3.1416,v=2.78,w=137.0,K=31.0])));
 
\qquad\qquad                                   .9999957122

\smallskip
 
$>$ simplify(eval(dif[1,3],[x=9,y=7,z=--5,

\qquad\qquad\qquad\qquad\qquad\ \ u=31416,v=278,w=137,K=31]));
 
\qquad\qquad                                  0
 
\smallskip

$>$ factor(simplify(eval(quot[1,3],[x=9,y=7,z=--5,

\qquad\qquad\qquad\qquad\qquad\qquad\qquad\ u=31416,v=278,w=137,K=31])));
 
\qquad\qquad                                  1

\smallskip
 
$>$ simplify(eval(dif[2,1],[x=131,y=17,z=--59,

\qquad\qquad\qquad\qquad\qquad\ \ u=61413,v=872,w=1/137,K=2]));
                                       
\qquad\qquad                                  0 

\smallskip

$>$ factor(simplify(eval(quot[2,1],[x=131,y=17,z=--59,

\qquad\qquad\qquad\qquad\qquad\qquad\qquad\ u=61413,v=872,w=1/137,K=2])));
 
\qquad\qquad                                  1

\smallskip
 
$>$ simplify(eval(dif[2,2],[x=1,y=2,z=--3,
 
\qquad\qquad\qquad\qquad\qquad\ \ u=71,v=5,w=1/137,K=29]));
                                        
\qquad\qquad                                  0
 
\smallskip

$>$ factor(simplify(eval(quot[2,2],[x=1,y=2,z=--3,

\qquad\qquad\qquad\qquad\qquad\qquad\qquad\ u=71,v=5,w=1/137,K=29])));
 
\qquad\qquad                                  1

\smallskip
  
$>$ simplify(eval(dif[2,3],[x=1.0,y=2.0,z=--3.0,

\qquad\qquad\qquad\qquad\qquad\ \ u=3.1416,v=2.78,w=137.0,K=29.0]));
 
\qquad\qquad                                  --.005050000000
  
\smallskip

$>$ factor(simplify(eval(quot[2,3],[x=1.0,y=2.0,z=--3.0,
                       
\qquad\qquad\qquad\qquad\qquad\qquad\qquad\ u=3.1416,v=2.78,w=137.0,K=29.0])));
 
\qquad\qquad                                  1.000000215

\smallskip
 
$>$ simplify(eval(dif[3,1],[x=1,y=2,z=3,

\qquad\qquad\qquad\qquad\qquad\ \ u=5,v=7,w=11,K=13]));
 
\qquad\qquad                                  0
 
\smallskip

$>$ factor(simplify(eval(quot[3,1],[x=1,y=2,z=3,

\qquad\qquad\qquad\qquad\qquad\qquad\qquad\ u=5,v=7,w=11,K=13])));
 
\qquad\qquad                                  1

\smallskip
  
$>$ simplify(eval(dif[3,2],[x=199.7,y=--2.4168,z=3.5,

\qquad\qquad\qquad\qquad\qquad\ \ u=59,v=79,w=119,K=357]));
 
\qquad\qquad                                    .00003253970835

\smallskip
 
$>$ factor(simplify(eval(quot[3,2],[x=199.7,y=--2.4168,z=3.5,
                               
\qquad\qquad\qquad\qquad\qquad\qquad\qquad\ u=59,v=79,w=119,K=357])));
 
\qquad\qquad                                    .9999963000

\smallskip
 
$>$ simplify(eval(dif[3,3],[x=1,y=2,z=--3,

\qquad\qquad\qquad\qquad\qquad\ \ u=71,v=5,w=1/137,K=29]));

\qquad\qquad                                   0
 
\smallskip

$>$ factor(simplify(eval(quot[3,3],[x=1,y=2,z=--3,

\qquad\qquad\qquad\qquad\qquad\qquad\qquad\ u=71,v=5,w=1/137,K=29])));
 
\qquad\qquad                                   1

As one can see, the numerical evidence is overwhelmingly in favor of our 
Randers metric having constant positive flag curvature $K$.

\end{document}